%% file: main_arxiv.tex
\documentclass[twocolumn,5p]{amsart}
\usepackage{amsaddr}
\usepackage{geometry}
\input{tex_files/preprocess_arxiv}
\usepackage{xcolor}
\usepackage{graphicx}
\graphicspath{
{figs/bojana}
}
\hypersetup{draft}

\newcommand{\mat}[1]{\mathchoice{\displaystyle\mathbf#1}
{\textstyle\mathbf#1}{\scriptstyle\mathbf#1}
{\scriptscriptstyle\mathbf#1}}

\title{Bayesian identification of a projection-based Reduced Order Model for Computational Fluid Dynamics}

\author{Giovanni Stabile\textsuperscript{1,*}}
\address{\textsuperscript{1}SISSA, International School for Advanced Studies, Mathematics Area, mathLab, via Bonomea 265 - 34136, Trieste, Italy}
\thanks{\textsuperscript{*}Corresponding Author.}
\email{gstabile@sissa.it} 
\author{Bojana Rosic\textsuperscript{2}}
\address{\textsuperscript{2}University of Twente, Applied Mechanics and Data Analysis, Drienerlolaan 5 - 7522 NB Enschede, Netherlands}
\email{b.rosic@utwente.nl}

\date{}

\dedicatory{}


\begin{document}

\twocolumn[
  \begin{@twocolumnfalse}
    \maketitle \thispagestyle{empty}
\begin{abstract} 
In this paper we propose a Bayesian 
method as a numerical way to correct and stabilise projection-based reduced order models (ROM) in computational fluid dynamics problems. The approach is of hybrid type, and consists of the classical proper orthogonal decomposition driven Galerkin projection of the laminar part of the governing equations, and Bayesian identification of the correction term mimicking both the turbulence model and possible ROM-related instabilities given the full order data. In this manner the classical ROM approach is translated to the parameter identification problem on a set of nonlinear ordinary differential equations. Computationally the inverse problem is solved with the help of the Gauss-Markov-Kalman smoother in both ensemble and square-root polynomial chaos expansion forms. To reduce the dimension of the posterior space, a novel global variance based sensitivity analysis is proposed.
\end{abstract}
  \end{@twocolumnfalse}
]

\input{tex_files/01_introduction}
\input{tex_files/02_mathMod}
\input{tex_files/03_ROM}

\input{tex_files/04_UQ}

\input{tex_files/05_numResults}

\input{tex_files/06_conclusions}

\section*{Acknowledgements} 
We acknowledge the support provided by the German Science Foundation (Deutsche Forschungsgemeinschaft, DFG) as part of priority programs SPP 1886 and SPP 1748.

\bibliographystyle{amsplain}
\bibliography{bib/bibfile} 
\end{document}

%% file: tex_files/01_introduction.tex
\section{Introduction}\label{sec:intro}

In last decades projection-based reduced order modelling demonstrated to be a viable way to reduce the computational
burden associated with the simulation of complex fluid dynamics problems
\cite{Hesthaven2016,Quarteroni2016,GALLETTI2004,Iollo1998}. The computational reduction is
particularly beneficial when a large number of simulations associated with different input
parameter values, are needed, or a reduced computational cost is required (real-time
control, uncertainty quantification, inverse problems, optimisation, \dots). 

The main issues of projection-based reduced order modelling are accuracy and stability
of the method. In particular, in fluid dynamics problems, which are the focus of this work, three different types of instability or inaccuracy issues are observed so far:
\begin{itemize}
\item lack of inf-sup stability of the resulting reduced basis spaces,
\item instability and inaccuracy associated with advection dominated problems,
\item and instability and inaccuracy associated with long-term time integrations.
\end{itemize}
The first type of instability, that is associated with the fact of considering both velocity and pressure terms at the reduced order level, can be resolved using a variety of different numerical techniques.
The supremizer stabilisation method \cite{Rozza2007,Ballarin2015,StaRo2018},
a stabilisation based on the use of a Poisson equation for pressure
\cite{StaHiMoLoRo2017,Akhtar2009} or Petrov-Galerkin projection strategies
\cite{Caiazzo2014598,Carlberg2017} stand out among others. Along these lines, 
it is a common practice to neglect the pressure contribution and to generate a reduced order model that accounts for the velocity
field only. For more details on reduced order models that also consider the pressure contribution the reader is
referred to \cite{StaRo2018} and references therein. 

The second type of instability is usually resolved by
adding extra stabilisation terms either of the Streamline Upwind Petrov-Galerkin
(SUPG) or Variational MultiScale
\cite{Bergmann2009Enablers,StaBaZuRo2019,Iliescu_vms,wang_turb} type. 
The long-term instability issues, on the other hand, can be fixed with the help of the constrained
projection \cite{Fick2018} or by introducing the spectral viscosity diffusion convolution operator \cite{Sirisup2004}.

The aim of this work is to extend the previously mentioned approaches in order to deal with the second and third type of inaccuracies and instabilities from a learning perspective. Following the least square idea presented in \cite{Mohebujjaman2018,Xie2018}, we propose to use a Bayesian strategy to identify correction terms associated
with the reduced differential operators derived from a standard POD-Galerkin approach with or without additional physical
constrains. A priori modelled as the tensor valued Gaussian random variable, the correction term is learned by exploiting the information contained in the partial observation, e.g. only a few time steps of the full order model response projected onto the reduced basis space. In this regard, the reduced order model (ROM) 
is substituted by the corresponding probabilistic one that accounts for the sparsity of both observations and correction terms. In this manner the ROM is not fitted to the observation, but learned. In order to reduce the learning time, the full Bayes' approach 
is replaced by an approximation of the posterior measure. As the interest lies only in the posterior mean, in this paper the learning is based solely on the estimation of the posterior conditional mean. Following \cite{Ro19}, the conditional expectation is 
an optimal projection of the random variable onto the space of all random variables consistent with the data, and therefore can be estimated directly. The estimation 
is however based on the approximation of the optimal projector, and can be achieved in a statistical manner by using the ensemble of particles \cite{Evensen:2006:DAE:1206873}, or by using the functional approximation approach \cite{Ro19}. As the latter one incorporates more information, the primary focus of this work is to exploit this method for the reduced modelling purpose. 

The paper is organised as follows, in \autoref{sec:fom} is described the general mathematical framework of the considered problem. \autoref{sec:rom} briefly recalls the projection based reduced order model with the additional correction terms, whereas in
\autoref{sec:bayesian} the approximate Bayesian inference techniques used to identify the correction terms are presented. Finally, in \autoref{sec:num_results} the methodology is applied on a fluid dynamics problem and in \autoref{sec:conclusions} conclusions and outlook are presented.

%% file: tex_files/02_mathMod.tex
\section{The mathematical model and the full order problem}\label{sec:fom}


Let be given the unsteady incompressible Navier-Stokes equations described in an Eulerian framework on a space-time domain $Q = \mathcal{G} \times [0,T] \subset \mathbb{R}^d\times\mathbb{R}^+, \textrm{ }d=2,3$ with the vectorial velocity field $\bm{u}:Q \to \mathbb{R}^d$ and the scalar pressure field $p:Q \to \mathbb{R}$ such that:
\begin{equation}\label{eq:navstokes}
\begin{cases}
\bm{u_t}+ \bm{\nabla} \cdot (\bm{u} \otimes \bm{u})- \bm{\nabla} \cdot 2 \nu \bm{\nabla^s} \bm{u}=-\bm{\nabla}p &\mbox{ in } Q,\\
\bm{\nabla} \cdot \bm{u}=\bm{0} &\mbox{ in } Q,\\
\bm{u} (t,x) = \bm{f}(\bm{x}) &\mbox{ on } \Gamma_{In} \times [0,T],\\
\bm{u} (t,x) = \bm{0} &\mbox{ on } \Gamma_{0} \times [0,T],\\ 
(\nu\nabla \bm{u} - p\bm{I})\bm{n} = \bm{0} &\mbox{ on } \Gamma_{Out} \times [0,T],\\ 
\bm{u}(0,\bm{x})=\bm{k}(\bm{x}) &\mbox{ in } T_0,\\            
\end{cases}
\end{equation}
hold. Here, $\Gamma = \Gamma_{In} \cup \Gamma_{0} \cup \Gamma_{Out}$ is the boundary of $\mathcal{G}$ and, is composed of three different parts $\Gamma_{In}$, $\Gamma_{Out}$ and $\Gamma_0$ that indicate, respectively, inlet boundary, outlet boundary and physical walls. The term $\bm{f}(\bm{x})$ depicts the stationary non-homogeneous boundary condition, whereas $\bm{k}(\bm{x})$ denotes the initial condition for the velocity at $t=0$.  

Furthermore, let the above system of equations be discretised
by the finite volume method with a segregated pressure-velocity coupling using the PIMPLE algorithm implemented in the OpenFOAM package \cite{OF}, further referred as the full order model (FOM). The PIMPLE algorithm is a pressure-velocity coupling approach for unsteady flows that merges the PISO \cite{Issa1986} and the SIMPLE \cite{Patankar1972} schemes, for more details please see \autoref{sec:num_results}.
After the discretised system of equations is solved, the velocity $\bm u (t)$ is collected at different time instants $t_i \in \{t_1, \dots, t_{N_t}\} \subset [0,T]$ and stored in the snapshots matrix:
\begin{equation}
\bm{\mathcal{U}} = [\bm{u}(t_1), \dots, \bm{u}(t_{N_t})] \in \mathbb{R}^{N_c \times N_t},
\end{equation}
in which $N_c$ denotes the number of spatial discretisation cells, and $N_t$ is the number of time instances. Furthermore, the previously defined matrix is used to construct the reduced model in a projection based manner as described in the following sections. However, note that the proposed methodology is general enough and does not depend on the previously chosen discretisation scheme.

%% file: tex_files/03_ROM.tex
\section{The projection based reduced order model}\label{sec:rom}
Assuming that the solution manifold of the full order discretised Navier-Stokes equations can be approximated by a reduced number of dominant modes (reduced basis space) $\bm{\mathcal{S}}_u = [\bm{\varphi}_1, \dots, \bm{\varphi}_{N_r}],$
one may employ a Proper Orthogonal Decomposition (POD) approach \cite{Lumley1967TheSO} to approximate the discretised solution of \autoref{eq:navstokes} as $$\bm u \approx \bm u_r = \sum_i^{N_r} a_i \bm \varphi_i.$$ By projecting the full order model onto $\bm{\mathcal{S}}_u$ in a Galerkin manner and neglecting the pressure contribution, one arrives at a lower dimensional nonlinear system of ordinary differential equations (ODEs):
\begin{equation}\label{eq:ODE}
\bm M_r \bm{\dot a} = \nu \bm A_r - \bm a^T \bm C_r \bm a,
\end{equation}
in which $\bm M_r$, $\bm A_r$, $\bm C_r$
represent the Gram matrix, the diffusion 2nd order tensor and the 3rd order convective tensor, respectively, all of which are computed as: 
\begin{equation}
\begin{split}\label{eq:red_matrices}
& (\bm M_r)_{ij} = \ltwonorm{\bm{\varphi_i},\bm{\varphi_j} } \mbox{, } (\bm{A_r})_{ij} = \ltwonorm{\bm{\varphi_i}, \bm{\nabla} \cdot 2\bm{\nabla^s\varphi_j}} \mbox{, } \\
& (\bm{C}_r)_{ijk} = \ltwonorm{\bm{\varphi_i}, \bm{\nabla} \cdot (\bm{\varphi_j} \otimes \bm{\varphi_k})}.
\end{split}
\end{equation}
Note that $\bm M_r$ is simplified due to the orthonormality of the POD modes as $\bm M_r=\bm I$. 

As already mentioned in \autoref{sec:intro}, the resulting system of ODEs in \autoref{eq:ODE} may lead to inaccurate results or instability issues. Moreover, in some cases, the straightforward Galerkin projection of the full order discretised differential operators onto the reduced basis space is not even possible. We refer for example to the case of commercial codes which only provide residuals and solution snapshots as the only available information. Therefore, one has to rely on an approximation of the reduced operators using, for example, an equivalent open source solver for the projection. Similarly, as stated in \cite{HiStaMoRo2019} the projection of the additional terms due to turbulence modelling may become inconvenient, and therefore only an approximation of the reduced order model can be made. The correction terms associated with the turbulence part of ROM can be then approximated with the help of a data-driven approach. In contrast to this, in a purely data-driven setting the structure of the reduced order model is usually directly deduced from the available data \cite{Benosman2017} or a ROM structure is first postulated, and then the reduced operators are identified by using data-driven methods \cite{Peherstorfer2015,Peherstorfer2016}. This can be challenging as the prior knowledge on such a structure is hard to define, whereas imposing the prior on the correction terms is more natural.
Following this, we add two correction terms to \autoref{eq:ODE} such that
\begin{equation}\label{eq:cor_ODE}
    \bm M_r \bm{\dot a} = \nu \bm A_r - \bm a^T \bm C_r \bm a + \nu \bm {\tilde{A}}_r - \bm a^T \bm{\tilde{C}}_r \bm a,
\end{equation}
holds. 
The term $\nu \bm {\tilde{A}}_r-\bm a^T \bm{\tilde{C}}_r \bm a$ aims to model all the neglected contributions and sources of error introduced by the Galerkin projection. As
these are generally not known, they are modelled as uncertain and further learned in a Bayesian framework given the measurement data, i.e. the time evolution of 
the FOM snapshots projected onto the POD modes:
\begin{equation}
    (\bm a_{FOM})_{ij}  = \ltwonorm{\bm{u}_{FOM}(t_i), \bm{\varphi}_j}.
\end{equation}
The final reduced order model is therefore a hybrid one that merges projection-based strategies with data driven techniques. 

%% file: tex_files/04_UQ.tex
\section{Bayesian identification of the correction terms}\label{sec:bayesian}

Given the nonlinear system of equations
\begin{equation}\label{eq:cor_ODE}
    \bm{\dot a} = \nu \bm A_r - \bm a^T \bm C_r \bm a + \nu \bm {\tilde{A}}_r - \bm a^T \bm{\tilde{C}}_r \bm a, 
\end{equation}
as described in the previous section, the goal is to find the correction parameters $\bm{q}$ 
represented by the elements of $\bm {\tilde{A}}_r$ and $\bm{\tilde{C}}_r$, such that the evolution 
of $ \bm{a}$ matches the evolution of the full order coefficients 
\begin{equation}\label{my_ode}
    (\bm{a}_{FOM})_{ij}, \quad i=1,..,N_t, j=1,..,N_r.
\end{equation}
In other words, we turn the projection based problem into the corresponding parameter estimation problem given the data $\mathbb{R}^{N_tN_r} \ni \bm y:=(\bm a_{FOM})_{ij}$. The estimation is in general an ill-posed problem as $\bm q$ is not directly observable, the nonlinear operator in \autoref{eq:cor_ODE} is not directly invertible, and the data
as well as the model in \autoref{eq:cor_ODE} are colored by an error. The sources of error are two-fold: the basis is generated by a truncated POD based approach which is accurate in a linear case only, and the full order model data are generated by an approximate discretisation technique. Therefore, the parameter estimation has to be regularised, and this is achieved here in a Bayesian setting. 
The prior information on $\bm q$---given in terms of probability density function $p_q(\bm q)$---acts as a regularisation term such that the problem of estimating the conditional probability density 
\begin{equation}\label{full_bayes}
    p_{\bm q|\bm y}(\bm q|\bm y)=\frac{p_{\bm y|\bm q}(\bm y|\bm q)p_q(\bm q)}{P(\bm y)}
\end{equation}
is well posed under the assumption that the joint probability density of $\bm q$ and $\bm y$ exists. Here, $p_{\bm y|\bm q}(\bm y|\bm q)$ denotes the likelihood function, whereas $P(\bm y)$ is the normalisation constant, or evidence. In other words, a priori  
the parameter $\bm q$ is modelled as a tensor-valued random variable in a probability space $L^2(\varOmega,\mathfrak{F},\mathbb{P};\mathbb{R}^s)$, $s={N_R}^2(1+{N_R})$ with the triple $\varOmega,\mathfrak{F},\mathbb{P}$ denoting the set of all events, the corresponding sigma-algebra, and the probability measure, respectively. For example, the 
diffusion correction term can be modelled as a tensor-valued Gaussian random variable, or a map
\begin{equation}
    \bm {\tilde{A}}_r(\omega): \varOmega \mapsto \mathbb{R}^{{N_R}\times {N_R}},
\end{equation}
and similarly $\bm{\tilde{C}}_r(\omega):\varOmega \mapsto  \mathbb{R}^{{N_R}\times {N_R}\times {N_R}}$. In general, the previous definition is not fully correct as 
the diffusion coefficient is known to be positive definite symmetric second order tensor. However, here we do not model the diffusion coefficient but its correction term (or error) which may equally take both positive and negative values, and therefore $\bm {\tilde{A}}_r(\omega)$ is not constrained on the positive definite manifold. Furthermore,  the data set $\bm y$ is assumed to be perturbed by a Gaussian noise (modelling error) with the zero mean and the covariance matrix $\bm C_\epsilon$, the value of which is assumed to be unknown 
and hence also estimated. 

As the parameter set $\bm q$ is assumed to be generated by $s$ independent random variables, the full Bayes' estimation as presented in \autoref{full_bayes} is not computationally tractable. To alleviate this, the posterior mean 
\begin{equation}\label{cond_exp}
    \mathbb{E}(\bm q | \bm y)=\int_\varOmega  \bm q p_{\bm q|\bm y}(\bm q|\bm y) \textrm{d} \bm q
\end{equation}
is further directly estimated
by using the projection-based algorithms as described in \cite{Ro19}. In a special case, when the noise and the prior distributions are both Gaussian, the estimate in \autoref{cond_exp} and in \autoref{full_bayes} are matching. 

Let $\phi: \mathbb{R}^s \mapsto \mathbb{R}$ be a strictly convex, differentiable function, and the corresponding loss function 
 $\mathcal{D}_\phi: \mathbb{R}^s \times \mathbb{R}
 \mapsto \mathbb{R}_+:=[0,+\infty)$ be defined as
 \begin{equation}\label{BLF}
  \mathcal{D}_\phi(\bm q,\hat{\bm q})=\mathcal{H}(\bm q)-\mathcal{H}(\hat{\bm q})=\phi(\bm q)-\phi(\hat{\bm q})-\langle \bm q-\hat{\bm q},\nabla \phi(\hat{\bm q})\rangle
 \end{equation}
 in which $\mathcal{H}(\bm q)=\phi(\hat{\bm q})+\langle \bm q-\hat{\bm q},\nabla \phi(\hat{\bm q}\rangle$ is the hyperplane tangent to $\phi$ at point $\hat{\bm q}$. The conditional expectation is then defined as the unique optimal projector for all loss functions \cite{Bregman1967}. 
\begin{equation}
\label{eq:optimality_blf}
 \bm q^*:=\mathbb{E}(\bm q|\mathfrak{B})=\underset{\hat{\bm q}\in L^2(\varOmega,\mathcal{B},\mathbb{P};\mathbb{R}^s)}{\arg \normalfont{\min}}
 \textrm{} \mathbb{E}(\mathcal{D}_\phi(\bm q,\hat{\bm q}))
\end{equation}
over all $\mathcal{B}$-measurable random variables $\hat{\bm q}$ in which $\mathcal{B}:=\sigma(\bm y)$ is the
sub-$\sigma$-algebra generated by a measurement $\bm y$.

Furthermore, assuming that $\phi$ takes the quadratic form, i.e. $\phi(\bm q)=\frac{1}{2}\|\bm q\|_{L^2}^2$, the loss $\mathcal{D}_\phi(\bm q,\hat{\bm q})$ in \autoref{eq:optimality_blf} reduces to the squared-Euclidean distance
\begin{equation}
 \mathcal{D}_\phi(\bm q,\hat{\bm q})=\|\bm q-\hat{\bm q}\|^2,
\end{equation}
and hence one may use the classical 
Pythagorean theorem to decompose the random variable $\bm q$ to the projected part 
$\bm q_p:=P_\mathcal{B}({\bm q})$, and its orthogonal component $\bm q_o:=(I-P_\mathcal{B}){\bm q}$. The projected component explains the observed data $\bm y$, whereas the residual information of our a priori knowledge remains in the orthogonal component, i.e.
\begin{equation}\label{ffa}
 {\bm q}_a=\mathbb{E}({\bm q}_f|{\bm y})+({\bm q}_f-\mathbb{E}({\bm q}_f|{\bm y}_f)).
\end{equation}
Thanks to this decomposition, one may form an update formula that can be further used to identify the ROM. Here, the indices $a$ and $f$ are used to denote the assimilated parameter and forecast (prior) parameter, respectively. The term $\bm y_f$ stands for the prediction of the observation obtained by solving \autoref{eq:cor_ODE} given a priori knowledge on the parameter $\bm q$ plus the modelling error $\bm \epsilon(\omega)$, i.e.
\begin{equation}\label{eq:cor_ODEn}
    \bm y_f(\omega)=\left(\bm a_{ij}(\bm q_f(\omega))+\bm \epsilon_{ij}(\omega)\right), \textrm{ }i=1,...,N_t, j=1,...,N_r.
\end{equation}
This is also the most expensive part of the update in \autoref{ffa} as the deterministic evolution for the coefficients $\bm a $ in \autoref{eq:cor_ODE}
turns into the stochastic one.

As the formula in \autoref{ffa} is relatively abstract, the conditional expectation is further materialised with the help of the measurable map $\varphi(\cdot): \bm y_f \mapsto \bm q_f$ according to the Doob-Dynkin's lemma \cite{bobrowski2005functional}. As the map is not known a priori, the easiest choice is to use the linear one due to computational simplicity. In other words, one assumes
\begin{equation}
\label{cond_map}
\mathbb{E}(\bm q_f|\bm y_f)\approx \bm K{\bm y_f}+{\bm b}
\end{equation}
in which the map coefficients $(\bm K,\bm b)$ are obtained by minimising the orthogonal component in \autoref{ffa}, i.e.
\begin{eqnarray}\label{inver_map_k}
 &&\underset{\bm K,\bm b}{\arg \min}\textrm{ } \mathbb{E}(\|\bm q_f-\mathbb{E}(\bm q_f|\bm y_f)\|_2^2)\nonumber\\
 &&=\underset{\bm K,\bm b}{\arg \min}\textrm{ } \mathbb{E}(\|\bm q_f-(\bm K\bm y_f+\bm b)\|_2^2).
\end{eqnarray}
From the optimality condition
\begin{equation}
 \forall \chi:\quad \mathbb{E}(\langle \bm q_f-(\bm K\bm y_f+\bm b),\chi \rangle) =0.
\end{equation}
 one obtains
the Kalman gain 
\begin{equation}
\label{kalmna_gain_eq}
\bm K={\bm C}_{{\bm q}_f,{\bm y}_f}({\bm C}_{{\bm y}_f})^{\dagger},
\end{equation}
and as a result of \autoref{ffa} a linear Gauss-Markov-Kalman (GMK) filter formula
\begin{equation}
\label{line_gaus_mar}
{\bm q}_a(\omega)={\bm q}_f(\omega)+{\bm K}({\bm y}-{\bm y}_f(\omega)).
\end{equation}
Here, $\dagger$ denotes the pseudo-inverse, ${\bm C}_{{\bm q}_f,{\bm y}_f}$ is the covariance
between the prior $\bm q_f$ and the observation forecast $\bm y_f$, and 
${\bm C}_{{\bm y}_f}={\bm C}_{\bm{a}(\bm q_f)}+\bm C_{\bm \epsilon}$ is the auto-covariance of $\bm y_f$ consisting of the
forecast covariance ${\bm C}_{\bm{a}(\bm q_f)}$ and the "measurement" covariance 
$\bm C_{\bm \epsilon}$, see \autoref{eq:cor_ODEn}. 

\subsection{Discretisation}

For computational purposes, \autoref{line_gaus_mar} is further discretised in both Monte Carlo sampling and functional approximation manner. These are known as the ensemble Kalman filter (EnKF) \cite{Evensen:2006:DAE:1206873} and the polynomial chaos expansion (PCE) GMK filter, see \cite{Ro19}. 
The EnkF updates the posterior variable by sampling 
\begin{equation}
\label{line_gaus_mar_samp}
{\bm q}_a(\omega_i)={\bm q}_f(\omega_i)+{\bm K}_s({\bm y}-{\bm y}_f(\omega_i)), \quad i=1,..,Z,
\end{equation}
in which $\bm K_s$ is estimated using 
\autoref{kalmna_gain_eq} and the formula for the statistical covariance
\begin{equation}
\label{covggg}
  \bm{C}_{{x,y}}=\frac{1}{N-1} \bm \tilde{X}\bm \tilde{Y}^T.
\end{equation}
Here, $\bm \tilde{X}$ is the vector of samples of the fluctuating part (without mean) of the random variable $x$. Similar holds for $\tilde{Y}$.

In a PCE update form, each of the random variables in \autoref{line_gaus_mar} is expressed in a coordinate system described by simpler standard random variables such as Gaussian ones, i.e.
\begin{equation}
\label{pcex}
\bm q_a(\omega)\approx \hat{\bm{q}}_a(\omega)=\sum_{\alpha \in \mathcal{J}_{\Psi}} \bm{q}_a^{(\alpha)}\Psi_\alpha(\bm{\xi}(\omega)),
\end{equation}
in which $\Psi_\alpha(\bm{\xi}(\omega))$ are the orthogonal polynomials and $\bm{\xi}:=(\xi_i)_{i=1}^m, m=s+N_tN_r$. Following this, 
\autoref{line_gaus_mar} rewrites to
\begin{equation}\label{line_gaus_mar1}
\bm q_a^{(\alpha)}={\bm q}_f^{(\alpha)}+\hat{\bm K}({\bm y}\otimes \bm e_\alpha-{\bm y}_f^{(\alpha)}).
\end{equation}
with $\bm e=(1,0,..,0)^T$ and $\hat{\bm K}$ being the disretised Kalman gain that is estimated using \autoref{kalmna_gain_eq} and 
\begin{equation}
\label{covggg}
  \bm {C}_{{x,y}}=\tilde{\bm {x}}\bm{\Delta}\tilde{\bm{y}}^T.
\end{equation}
 Here, 
$\tilde{\bm{x}}:=\{...,\bm{x}^{(\alpha)},...\}^T, \alpha>0$ and 
$\bm{\Delta}:=\mathbb{E}(\Psi_\alpha \Psi_\beta)$.

However, the previous approximation can be high-dimensional as it 
incorporates the prediction of the measurement error, see \autoref{eq:cor_ODEn}. If $\bm \epsilon$ is described by a sequence of
i.i.d.~Gaussian random variables,
then the posterior random variable lives in a high-dimensional probability space. To avoid this type of problem, the posterior random variable is further constrained to 
live in the prior approximation space. In other words the goal is to map the prior random variable approximated by a PCE $\hat{\bm q}_f$
to the posterior random variable approximated by $\hat{\bm q}_a$ such that the cardinality of both PCEs is equal. 
As the PCEs  $\hat{\bm q}_f,\hat{\bm q}_a$ directly relate to the covariance matrices,
the previous goal is same as finding an optimal linear map between the square-root\footnote{We use 
term square root for $S_f$, even though this is not the square-root of the matrix $ \bm{C}_{{q_f}}$ but its factorisation.} of the prior covariance matrix 
\begin{equation}
 \bm {S}_f=\tilde{\bm q}_f \bm{\Delta}^{-1/2},
\end{equation}
and the posterior one 
\begin{equation}
 \bm {S}_a=\tilde{\bm q}_a \bm{\Delta}^{-1/2},
\end{equation}
such that
\begin{equation}
\label{rel_cov}
 \bm {S}_a=\bm {S}_f\bm {W}^T
\end{equation}
holds.
Here, $\bm {W}$ denotes the corresponding transformation matrix that is not known, and can be found directly from the updating formula, for more details please see \cite{pajonk_bojana}.

Once the posterior square-root covariance function is evaluated, one may estimate the posterior fluctuating part
according to
\begin{equation}
\label{pce_exp}
\tilde{\bm{q}}_a=\bm{S}_a\bm{\Delta}^{-1/2}=\bm{S}_f(\bm{W})^T\bm{\Delta}^{-1/2},
\end{equation}
which together with the mean part
\begin{equation}
 \bar{\bm{q}}_{a}=\bar{\bm{q}}_{f}+\hat{\bm{K}}(\bm{y}-
 \bar{\bm{y}}_{f})
\end{equation}
builds the complete posterior PCE. 

\subsection{Uncertainty quantifcation via PCE}\label{subsec:uq}

The estimation of the reduced order model greatly depends on the measurement prediction $\bm y_f$, as well as its PCE proxy form $\hat{\bm y}_f$. 
To reduce the overall computational burden of the proxy,
 the propagation of the uncertainty through the forward problem
 can be achieved in 
a data-driven non-intrusive setting. This is here carried out by assuming the ansatz for the predicted measurement in a form 
\begin{equation}
\label{regression}
 \hat{\bm{y}}_f(\omega)=\sum_{\alpha \in \mathcal{J}} \bm{y}_f^{(\alpha)}\varPsi_\alpha(\bm{\xi}(\omega))=\mat{\Psi}\bm{v},
\end{equation}
and estimating the unknown coefficients 
$\bm{v}$ given $N$ samples $\bm w:=(\bm{y}_f(\omega_i))_{i=1}^N$, i.e.~
\begin{equation}\label{pce}
 \bm{y}_f(\omega_i)=\sum_{\alpha \in \mathcal{J}} \bm{y}_f^{(\alpha)}\varPsi_\alpha(\bm{\xi}(\omega_i))
\end{equation}
for $i=1,...,N$. In other words, the goal is to solve the system
\begin{equation}
 \label{reg}
 \bm{w}=\mat{\Psi} \bm{v}
\end{equation}
in which the number of samples $N\leq P:=\textrm{card } \mathcal{J}$, the cardinality of the PCE. As the system in \autoref{reg}
is undetermined, the solution can be only obtained if additional information is available. Having that a priori knowledge on the 
current observation exists, one may model 
the unknown coefficients $\bm{v}$ a priori as random variables
in $L^2(\varOmega_v,\mathcal{F}_v,\mathbb{P}_v;\mathbb{R}^P)$, i.e.
\begin{center}
 $\bm{v}(\omega_v): \varOmega_v \rightarrow \mathbb{R}^P,$
\end{center}
and further use the Laplace prior
\begin{center}
$ \bm{v} \sim e^{-{\|\bm{v}\|_1}} $
\end{center}
to model the coefficients and promote their sparsity. As the work with a Laplace prior is computationally difficult, in this paper the Laplace distribution is simulated by 
using the corresponding hyperprior as advocated in a relevance vector machine approach, the method used in this paper, \cite{tipping}. 
The hyperprior is modelled as
$$p(\bm v|\varpi)=\prod_{\alpha \in \mathcal{J}} p(\bm v^{(\alpha)}|\varpi_\alpha), \quad \bm v^{(\alpha)}\sim \mathcal{N}(0,\varpi_\alpha^{-1})$$
with $\varpi_\alpha$ being the precision of each PCE coefficient modelled by Gamma prior $p(\varpi_\alpha)$.
By marginalising over $\varpi$ one obtains the overall prior
 $$p(\bm v)=\prod_{\alpha \in \mathcal{J}} \int   p(\bm v^{(\alpha)}|\varpi_\alpha) p(\varpi_\alpha) \textrm{d}\varpi_\alpha$$
 which is further simplified by taking the most probable values for $\bm {\varpi}$, i.e.~$\bm {\varpi}_{MP}$. To estimate the coefficients in \autoref{reg} we further use the Bayes's rule in a form 
\begin{equation}\label{bayessparse}
p(\bm v,\bm \varpi,\bm \sigma^2|\bm {w})=\frac{p(\bm {w}|\bm v,\bm {\varpi},\bm {\sigma}^2)}{p(\bm {w})} p(\bm v,\bm {\varpi},\bm {\sigma}^2)
\end{equation}
 in which the coefficients $\bm v$, the precision $\bm {\varpi}$ and the regression error\footnote{This is the error of the proxy model, here assumed to be Gaussian, i.e. $\mathcal{N}(0,\sigma^2).$} $\bm {\sigma}^2$ are assumed to be uncertain.
For computational reasons the posterior is further factorised into
$$p(\bm v,\bm {\varpi},\bm {\sigma}^2|\bm {w})=p(\bm v|\bm {u},\bm {\varpi},\bm {\sigma}^2)p(\bm {\varpi},\bm {\sigma}^2|\bm {w})$$
in which the first factoring term is the convolution of normals
$p(\bm v|\bm {w},\bm {\varpi},\bm {\sigma}^2) \sim \mathcal{N}(\bm {\mu},\bm {\Sigma})$, whereas
the second factor $p(\bm {\varpi},\bm {\sigma}^2|\bm {w})$ cannot be computed analytically,
and thus is approximated by delta function $p(\bm {\varpi},\bm {\sigma}^2|\bm {u})\approx \delta (\bm {\varpi}_{MP},\bm {\sigma}_{MP})$.
The estimate $(\bm {\varpi}_{MP},\bm {\sigma}_{MP})$ is obtained from
$$p(\bm {\varpi},\bm {\sigma}^2|\bm {w})\propto p(\bm {u}|\bm {\varpi},\bm {\sigma}^2)p(\bm {\varpi})p(\bm {\sigma}^2)$$
by maximising the evidence (marginal likelihood)
$$p(\bm {w}|\bm {\varpi},\bm {\sigma}^2)=\int p(\bm {w}|\bm v,\bm {\sigma}^2)p(\bm v|\bm {\varpi}) \textrm{d}\bm v$$
as further discussed in detail in \cite{tipping}.

\subsection{Sensitivity analysis}\label{subsec:sensitivity}
\begin{figure}
\includegraphics[width=0.45\textwidth]{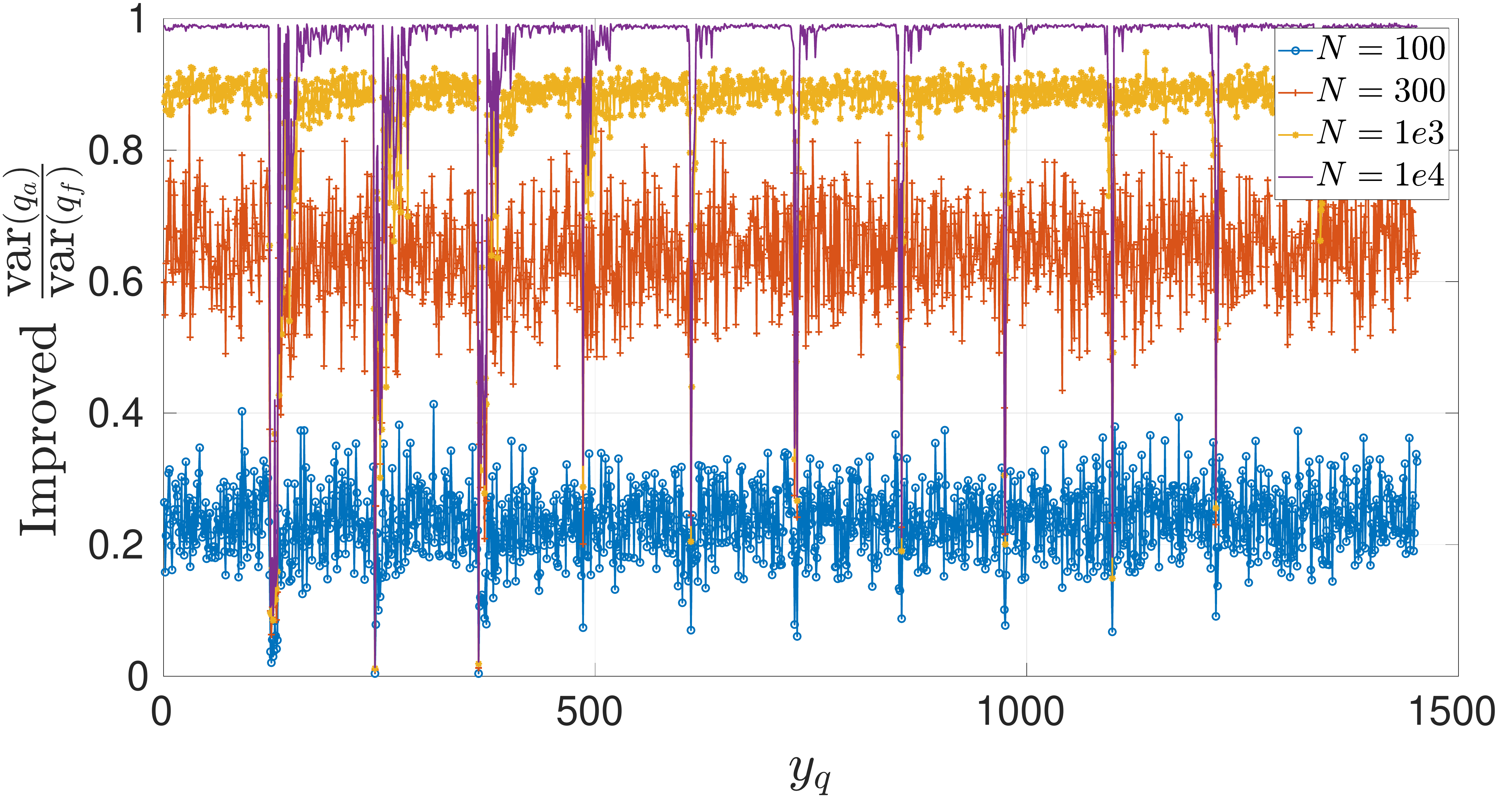}
\caption{Sensitivity analysis using classical EnKF}\label{fig:fig1imp}
\end{figure}

\begin{figure}
\includegraphics[width=0.45\textwidth]{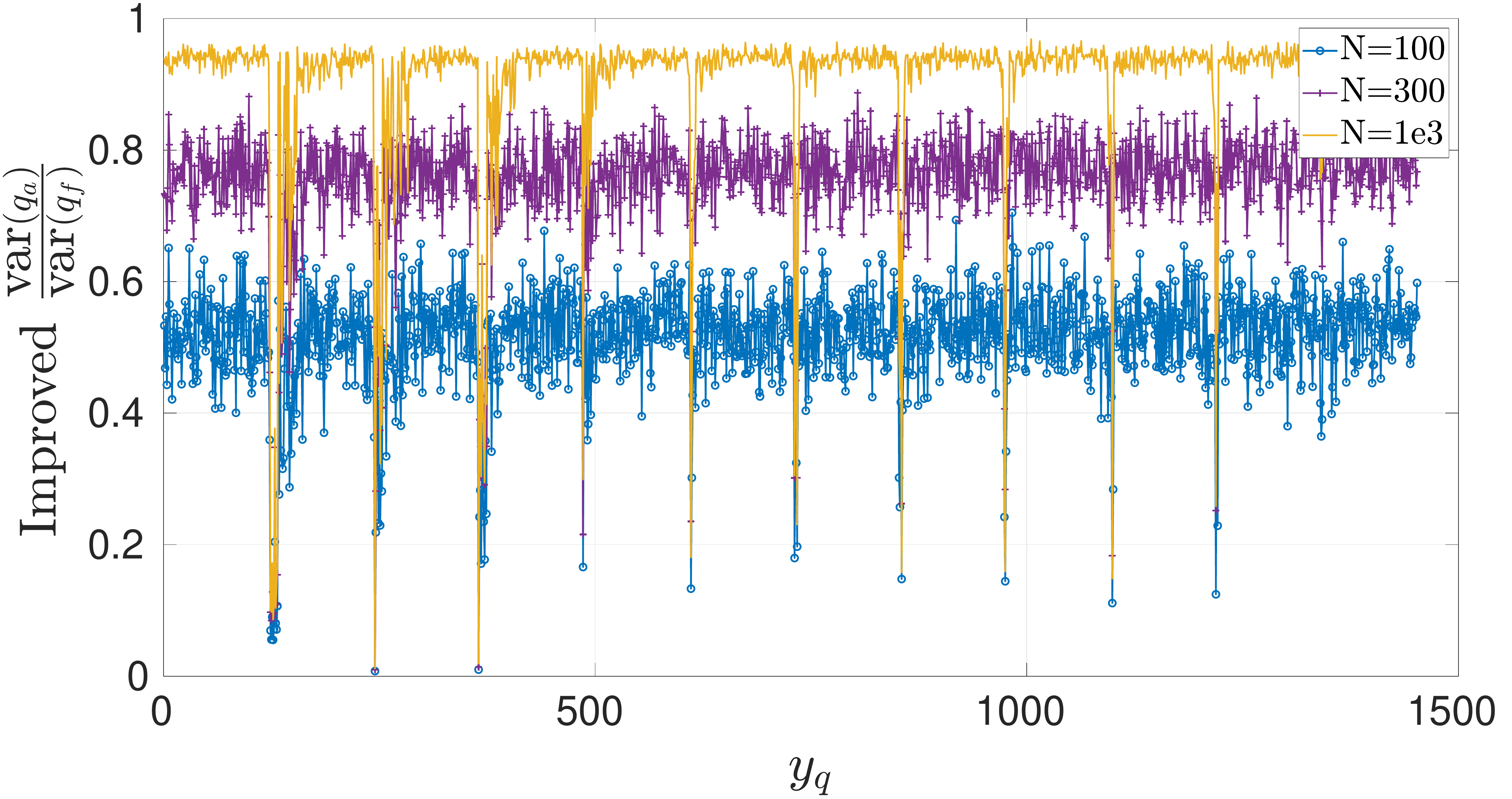}
\caption{Sensitivity analysis using improved EnKF}\label{fig:fig2imp}
\end{figure}

Even though the posterior variable is constrained to lie on the lower dimensional prior probability space, the PCE in \autoref{pce_exp} is not low-dimensional per se. The reason is that the parameter $\bm q$ is modelled as a tensor valued-random variable made of independent Gaussian variables. Due to independency relation, the PCE has to account for $s$ random variables $\xi_i, i=1,...,s$. This may lead to potentially high PCE cardinality that is not computationally feasible. To further reduce the computational burden, one may perform sensitivity analysis and discard variables $\xi_j,j=1,..\ell$ that are not contributing to the posterior. One way of achieving this is to use the variance based techniques. In general practice, the variance based 
sensitivity analysis a la Sobol \cite{sobol} requires estimation of the variance of the predicted measurement $\textrm{var } \bm y_f$ given all random variables of consideration $\bm \xi$, as well as the effect of the particular random variable $\xi_i, i=1,...,s$. The dimension reduction is then 
judged given the first Sobol's index defined as 
\begin{equation}
    S_i^\xi:=\frac{\textrm{var} (\mathbb{E}(\bm y_f|\xi_{i}))}{\textrm{var }\bm y_f}.
\end{equation}
As $q_{f,i}$ is only linear transformation of the random variable $\xi_i$, one may further rewrite
the previous equation as
\begin{equation}\label{sobol_q_f}
    S_i^{q_f}:=\frac{\textrm{var} (\mathbb{E}(\bm y_f|q_{f,i}))}{\textrm{var } \bm y_f}.
\end{equation}
If $S_i^{q_f}$ is low, then one may discard random variable $q_{f,i}$. In other words, $\bm y_f$ is not sensitive on $q_{f,i}$. However, this would also mean that $q_{f,i}$
is not sensitive on $\bm y_f$. In other words, 
one may also look at the dual definition of the previous equation
\begin{equation}
    J^{y_f}:=\frac{\textrm{var} (\mathbb{E}(\bm q_{f}| \bm y_{f,(k)}))}{\textrm{var }  \bm q_{f}}
\end{equation}
and estimate sensitivity of $\bm q_{f}$ given the forecast $\bm y_{f,(k)}:=\bm a_{(:,k)}+\bm \epsilon_{(:,k)}, k=1,..,N_r$. Following \autoref{cond_map}, one may further write 
\begin{equation}
    J^{y_f}:=\frac{\textrm{var} (\bm K_k \bm y_{f,(k)}+ \bm b| \bm y_{f,(k)})}{\textrm{var } \bm q_{f}},
\end{equation}
and conclude that if $\bm K_k$ is small enough, $\bm q_{f}$ can be considered as not sensitive on $\bm y_{f,(k)}$, and vice versa.
This would further mean that the posterior $\bm q_{a,k}$ would not be updated according to \autoref{line_gaus_mar}. Hence, one may use the update in \autoref{line_gaus_mar_samp} to perform sensitivity analysis according to the modified index
\begin{equation}\label{imp_ratio}
    J_k^{q_a}:=\frac{\textrm{var} (q_{a,k})}{\textrm{var }  q_{f,k}}.
\end{equation}
If the index is high, the posterior is close to the prior, and vice versa. 
However, as in any statistical approach, the EnKF estimate of $J_k^{q_a}$ is not reliable unless in the limit, see \autoref{fig:fig1imp}. 
The prior and forecast samples are correct as they are directly estimated from their corresponding models. However, the Kalman gain in \autoref{line_gaus_mar_samp}
is the only one prone to the statistical errors due to the approximation of the covariance matrices, see \autoref{covggg}. The error appears in the forecast covariance matrix $\bm C_{\bm y_f}$ and the cross-one $\bm C_{\bm q_f,\bm y_f}$, whereas the prior covariance $\bm C_{\bm q_f}$ matrix is known precisely. The former two can be estimated given the mapping $\bm H: \bm q_f \mapsto \bm y_f$ between the prior and the forecast, which is further estimated in a similar manner as the PCE coefficients, see \autoref{reg}, by using the sparse Bayesian learning.
Once this is achieved, one may estimate the Kalman gain as:
\begin{equation}
   \bm K_{imp}:=\bm H \bm C_{\bm q_f} (\bm H\bm C_{\bm q_f} \bm H^T+\bm C_\epsilon)^\dagger,
\end{equation}
and hence obtain improved EnKF estimate and sensitivity analysis, see \autoref{fig:fig2imp}. Note that here we only observed the first order indices, even though the second order ones could also play a role, but this is the topic of another paper.

%% file: tex_files/05_numResults.tex
\section{Numerical Results}\label{sec:num_results}
\begin{figure}
\includegraphics[width=0.48\textwidth]{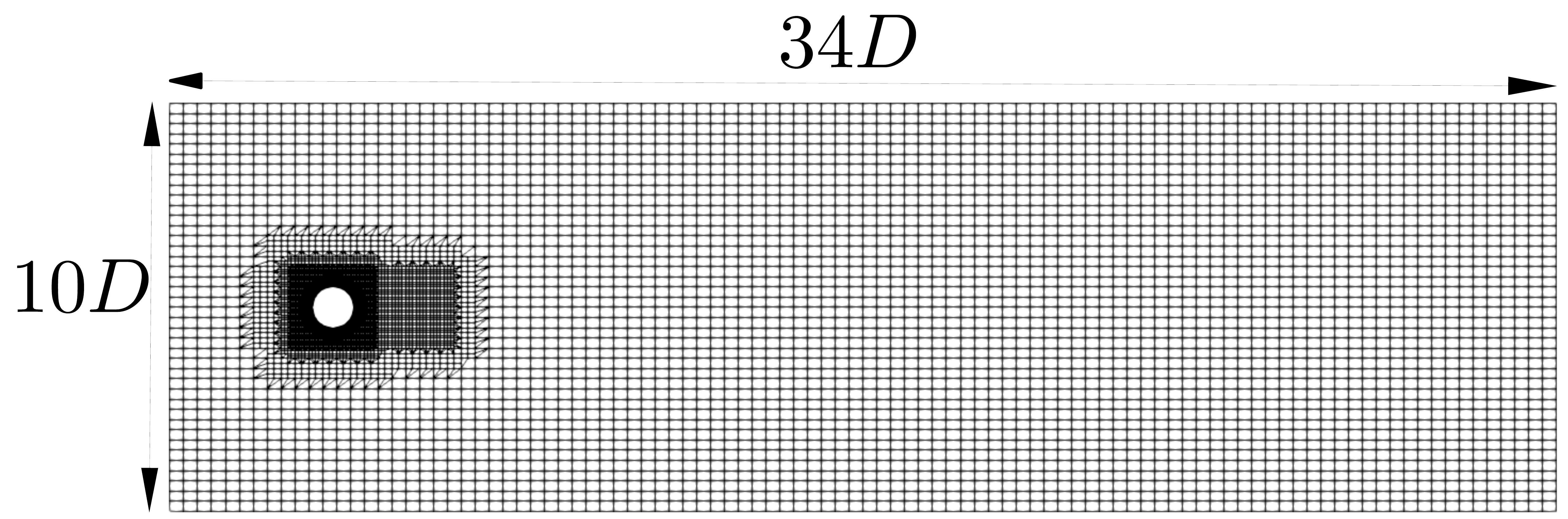}
\caption{Sketch of the computational mesh used for the full order simulations. The cylinder diameter is $D=1\si{m}$.}\label{fig:com_mesh}
\end{figure}

The proposed methodology is tested on the two-dimensional flow around a circular cylinder immersed in a rectangular domain, the dimensions of which are depicted in \autoref{fig:com_mesh}.
At the full order level a finite volume approximation counting $11644$ cells with a $k$-omega turbulence model is used. The kinematic viscosity $\nu=0.005 \si{m^2/s}$ considering an inlet velocity of $\bm u=1 \si{m/s}$ makes the Reynolds number equal to $Re=200$. At the inlet a uniform constant velocity is enforced while at the cylinder surface a no-slip condition is applied. The initial condition for both velocity and pressure are taken as homogeneous.

\begin{figure}
\centering
\begin{minipage}[c]{0.23\textwidth}
\includegraphics[width=\textwidth]{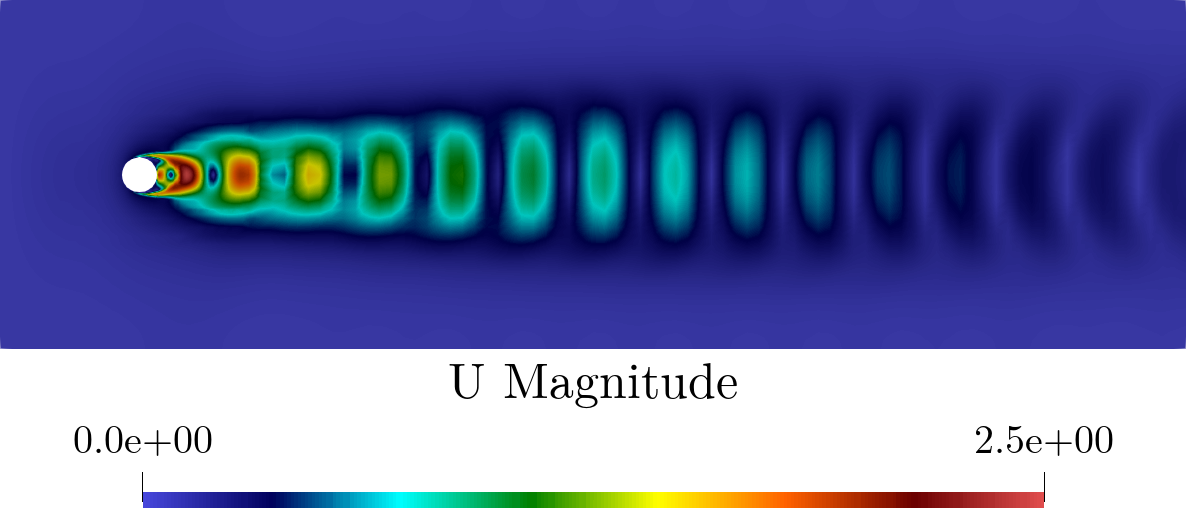}
\end{minipage}
\begin{minipage}[c]{0.23\textwidth}
\includegraphics[width=\textwidth]{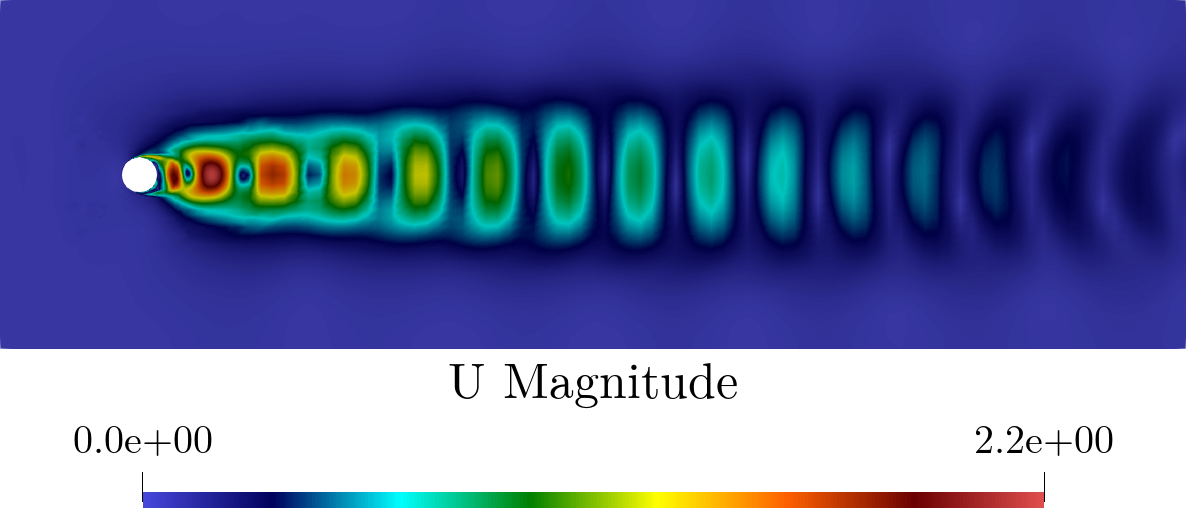}
\end{minipage}
\begin{minipage}[c]{0.23\textwidth}
\includegraphics[width=\textwidth]{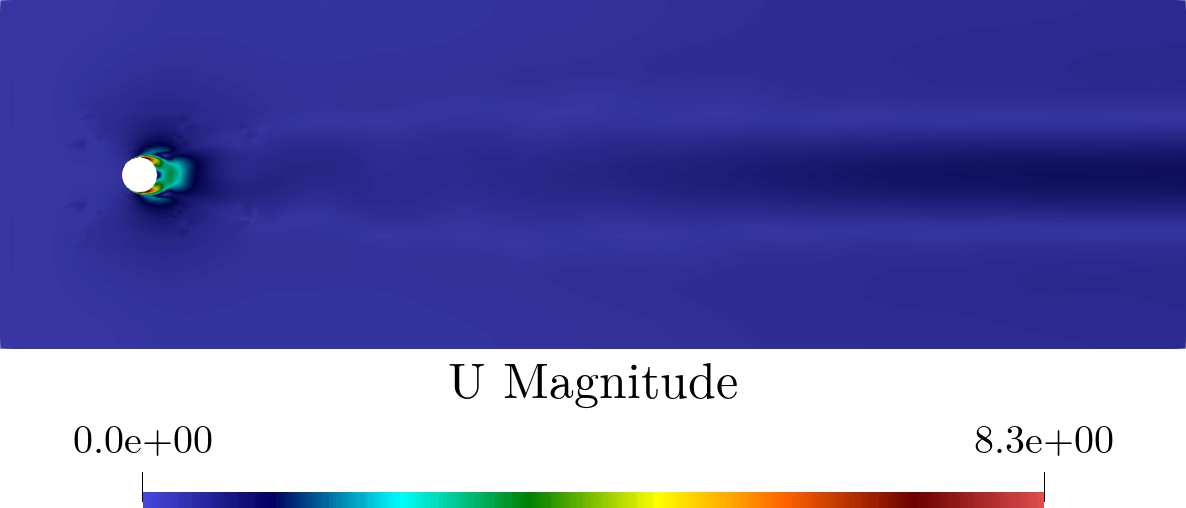}
\end{minipage}
\begin{minipage}[c]{0.23\textwidth}
\includegraphics[width=\textwidth]{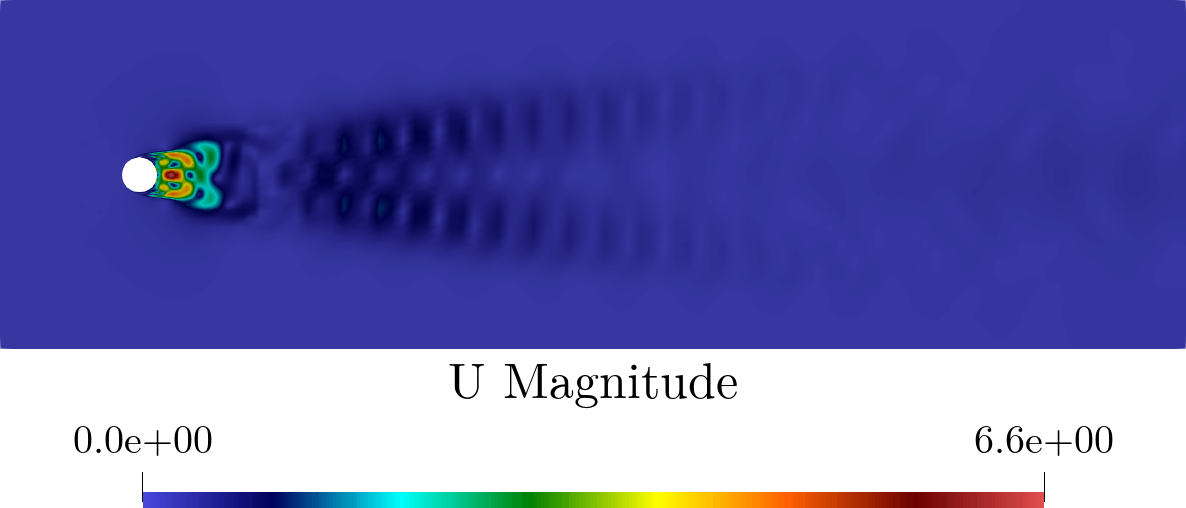}
\end{minipage}\caption{The first $4$ POD modes used to approximate the velocity field.}\label{fig:modes}
\end{figure}

\begin{figure*}
\centering
\begin{minipage}[c]{0.45\textwidth}
\includegraphics[width=\textwidth]{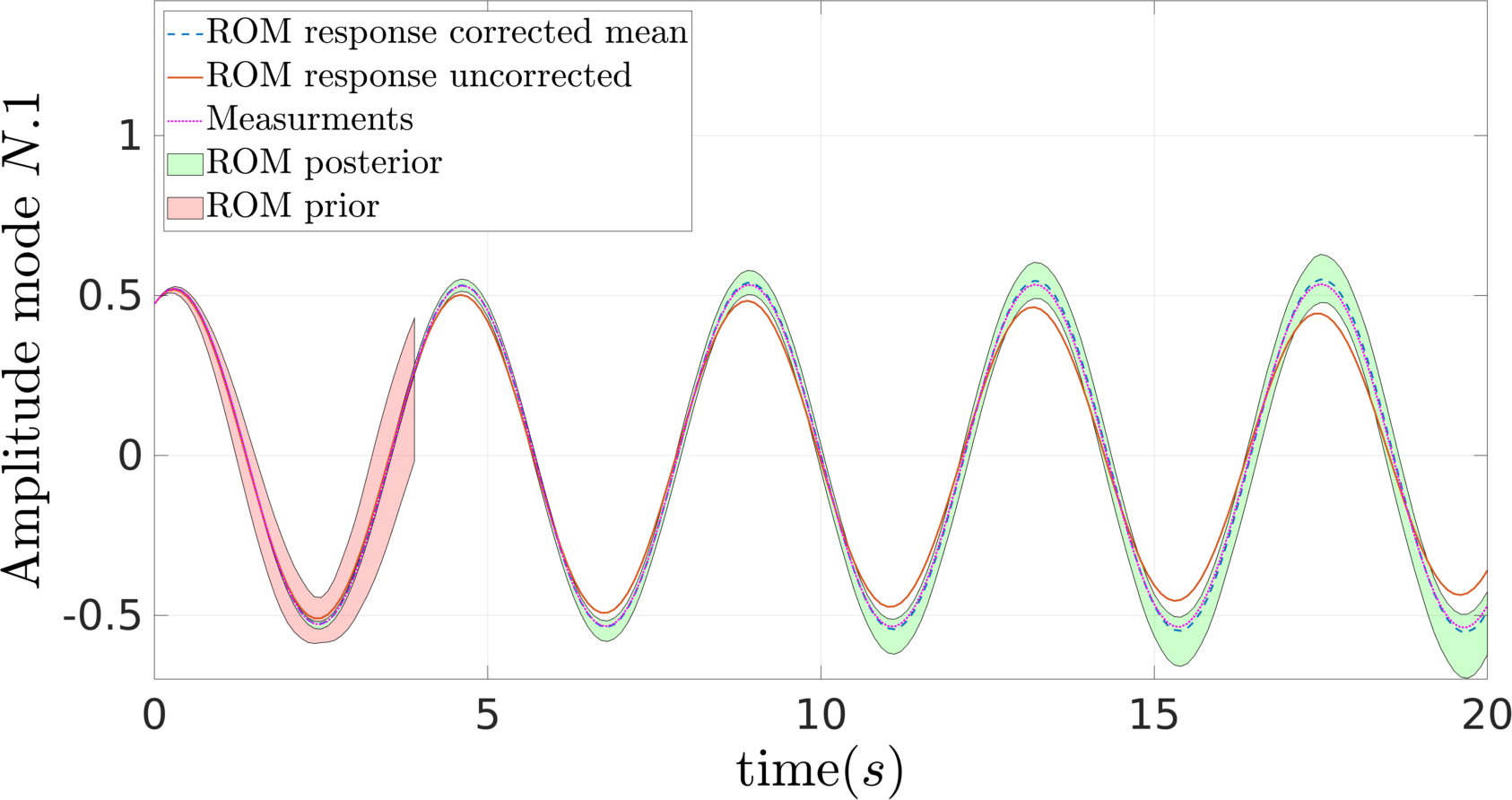}
\end{minipage}
\begin{minipage}[c]{0.45\textwidth}
\includegraphics[width=\textwidth]{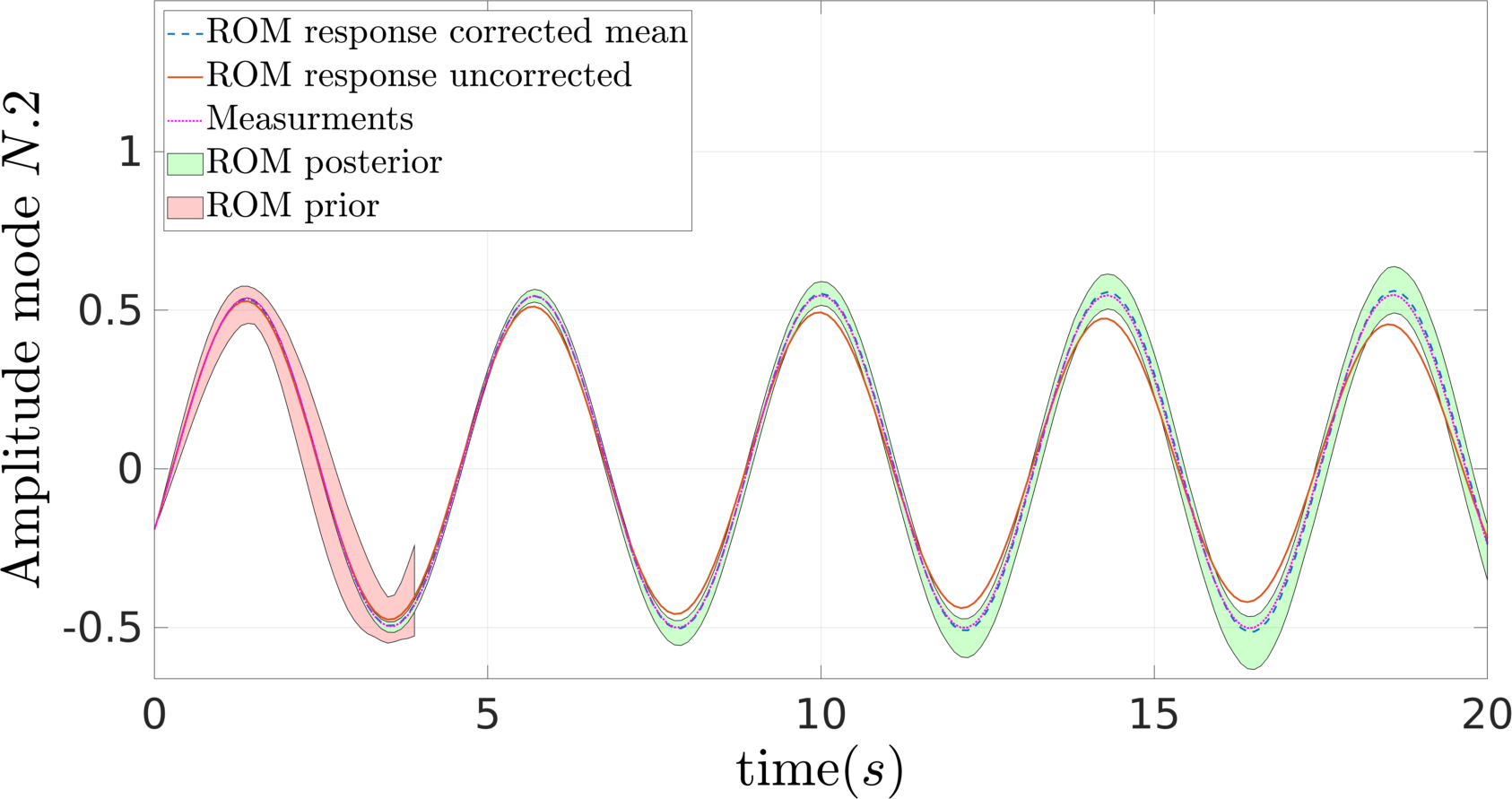}\end{minipage}\\
\begin{minipage}[c]{0.45\textwidth}
\includegraphics[width=\textwidth]{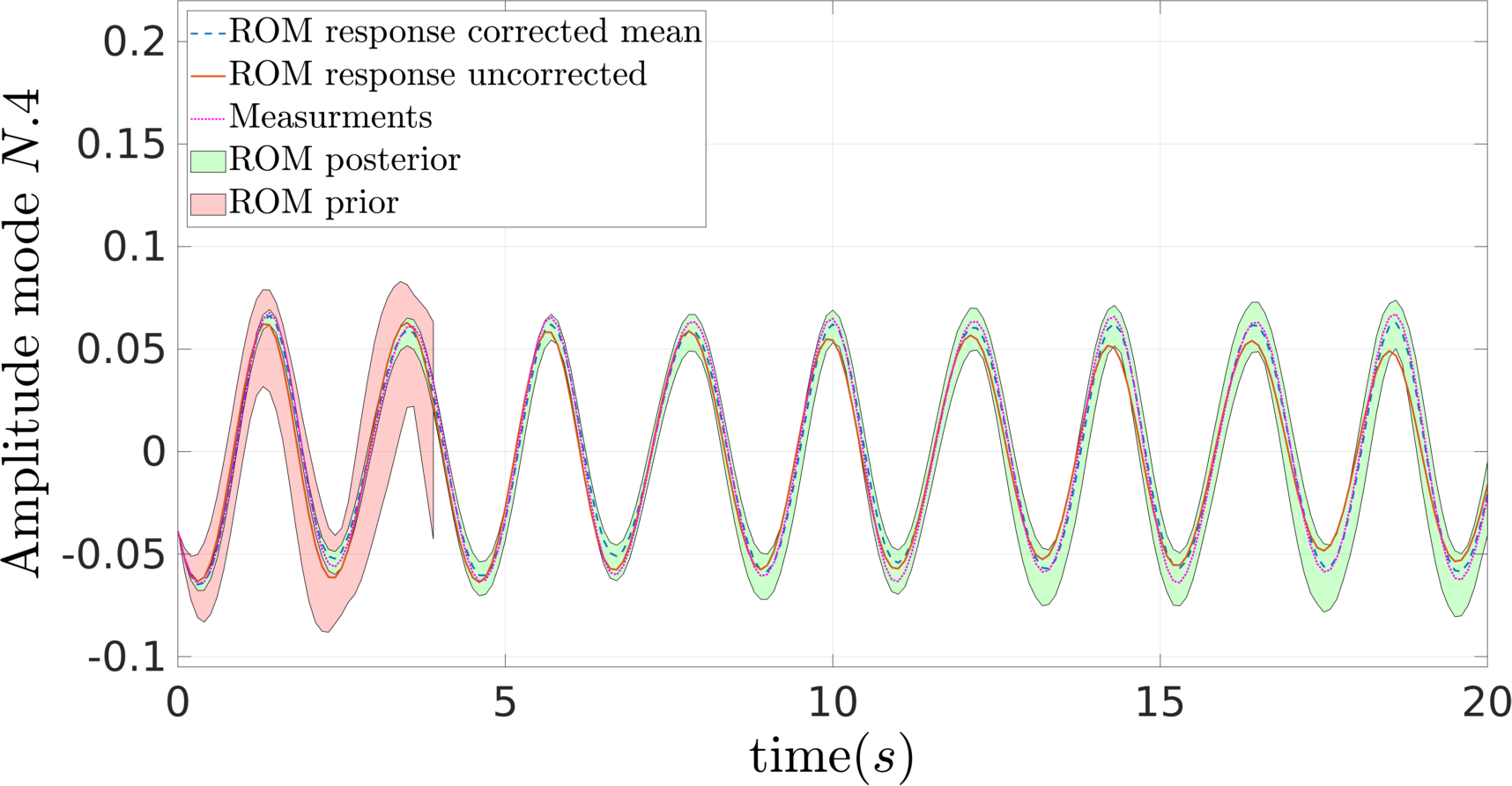}\end{minipage}
\begin{minipage}[c]{0.45\textwidth}
\includegraphics[width=\textwidth]{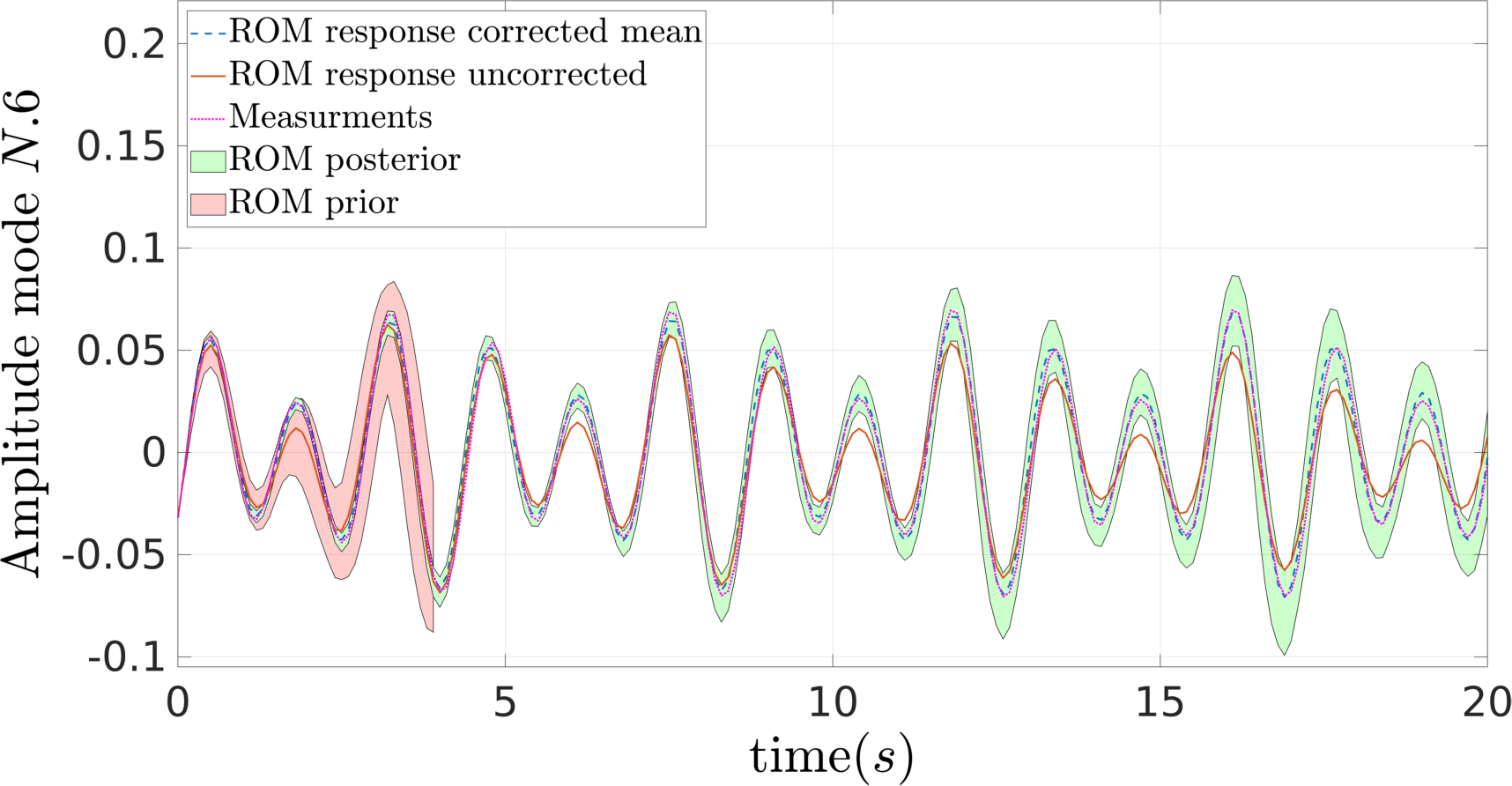}\end{minipage}
\caption{Probabilistic evolution of ROM coefficients for the modes $1,2,4$ and $6$. The plots depicts the $99\%$ quantiles of the model response before and after the identification, the uncorrected ROM (solid red line), the $L^2$-projection of the snapshots onto the POD modes (dotted purple line), and the mean posterior model response (dashed blue line).}\label{fig:coeefs}
\end{figure*}

\begin{figure}
    \centering
    \includegraphics[width=0.45 \textwidth]{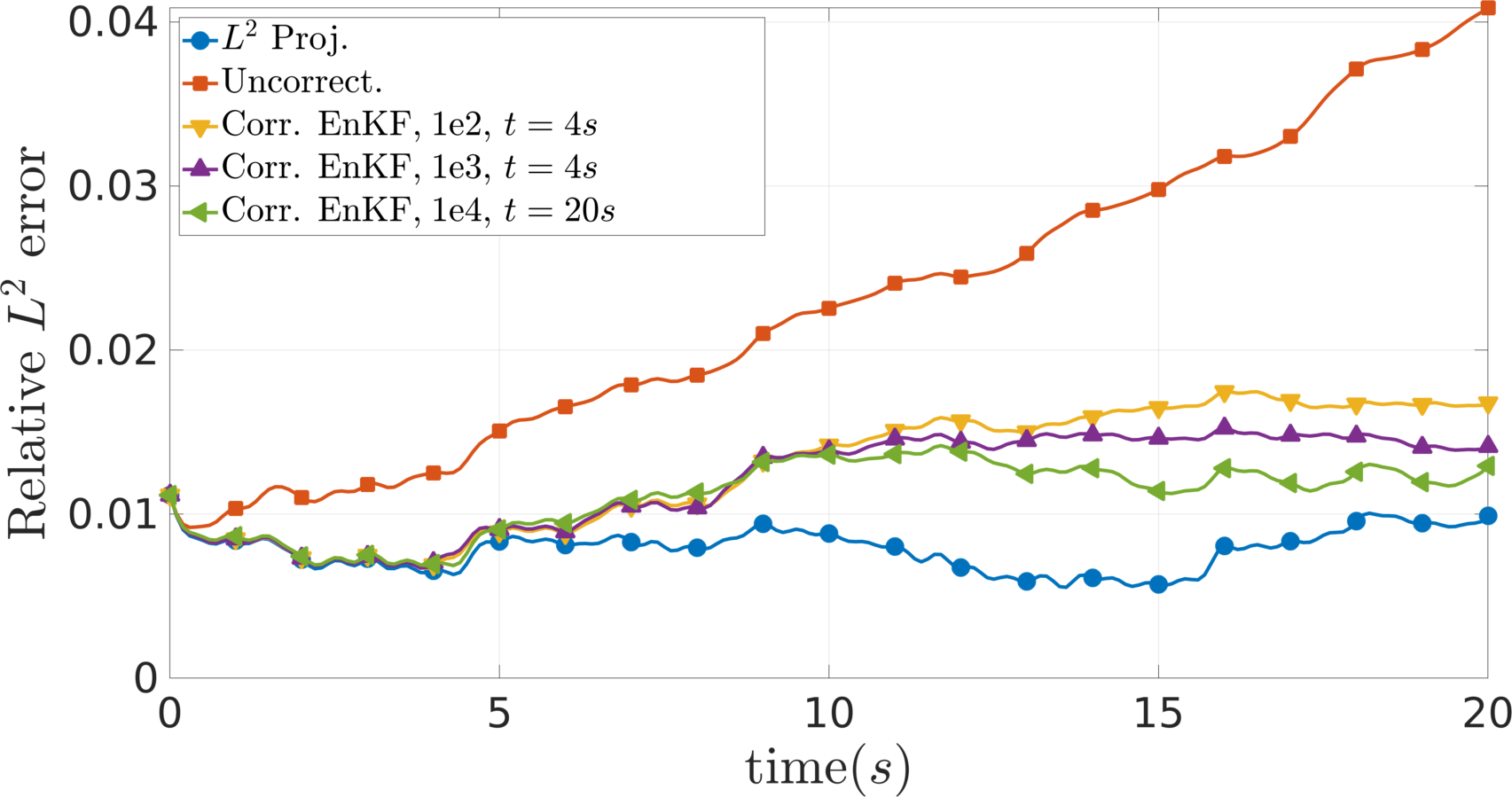}
    \caption{Time evolution of the relative error $\epsilon_u(t)$ of the Bayesian ROM using the EnKF smoother.}
    \label{fig:L2EnKF_all}
\end{figure}

\begin{figure}
    \centering
    \includegraphics[width=0.45 \textwidth]{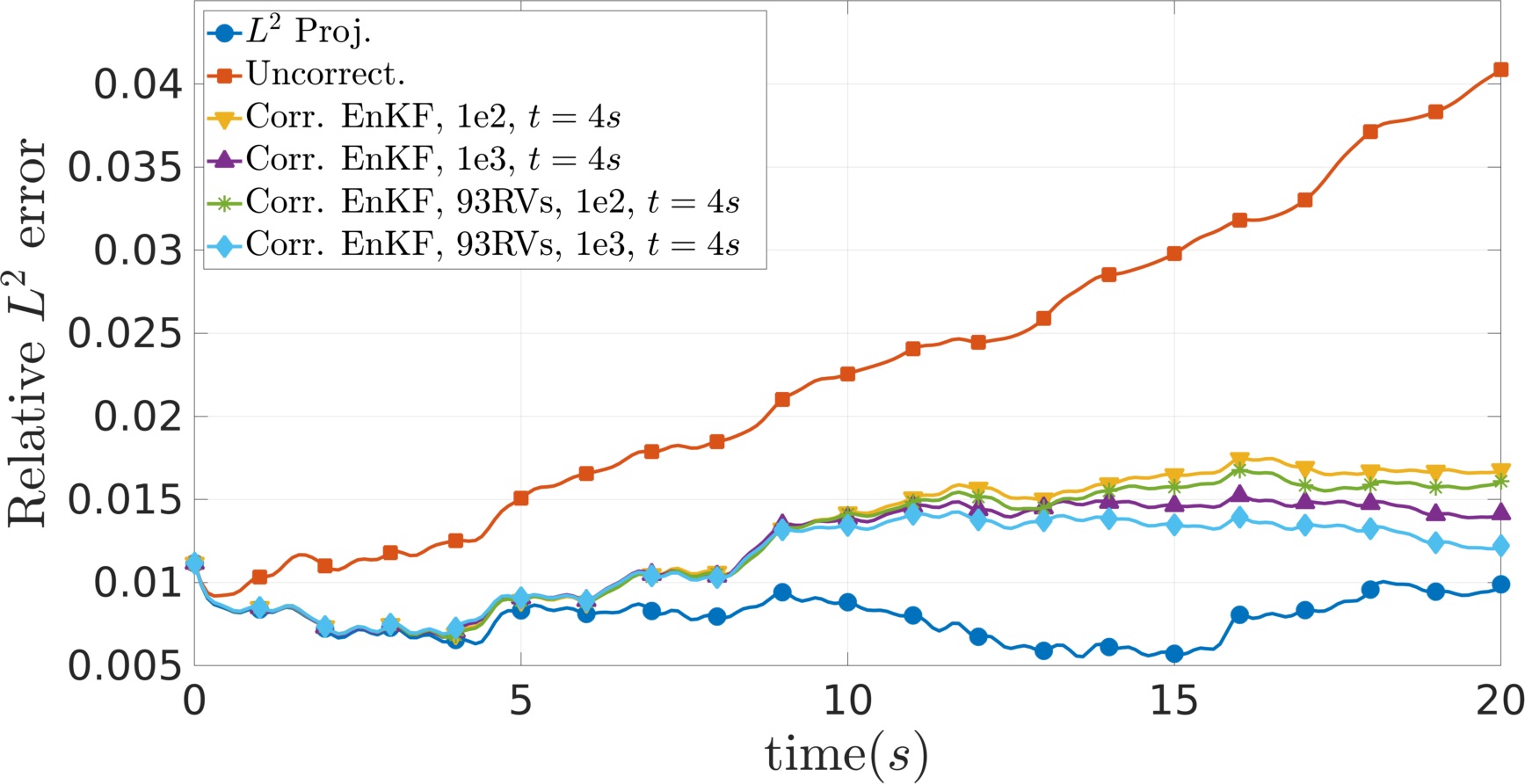}
    \caption{Time evolution of the relative error $\epsilon_u(t)$ of teh Bayesian ROM using the EnKF update technique on a reduced probabilistic space including $93$ random variables.}
    \label{fig:L2EnKF_sens}
\end{figure}

\begin{figure}
    \centering
    \includegraphics[width=0.45 \textwidth]{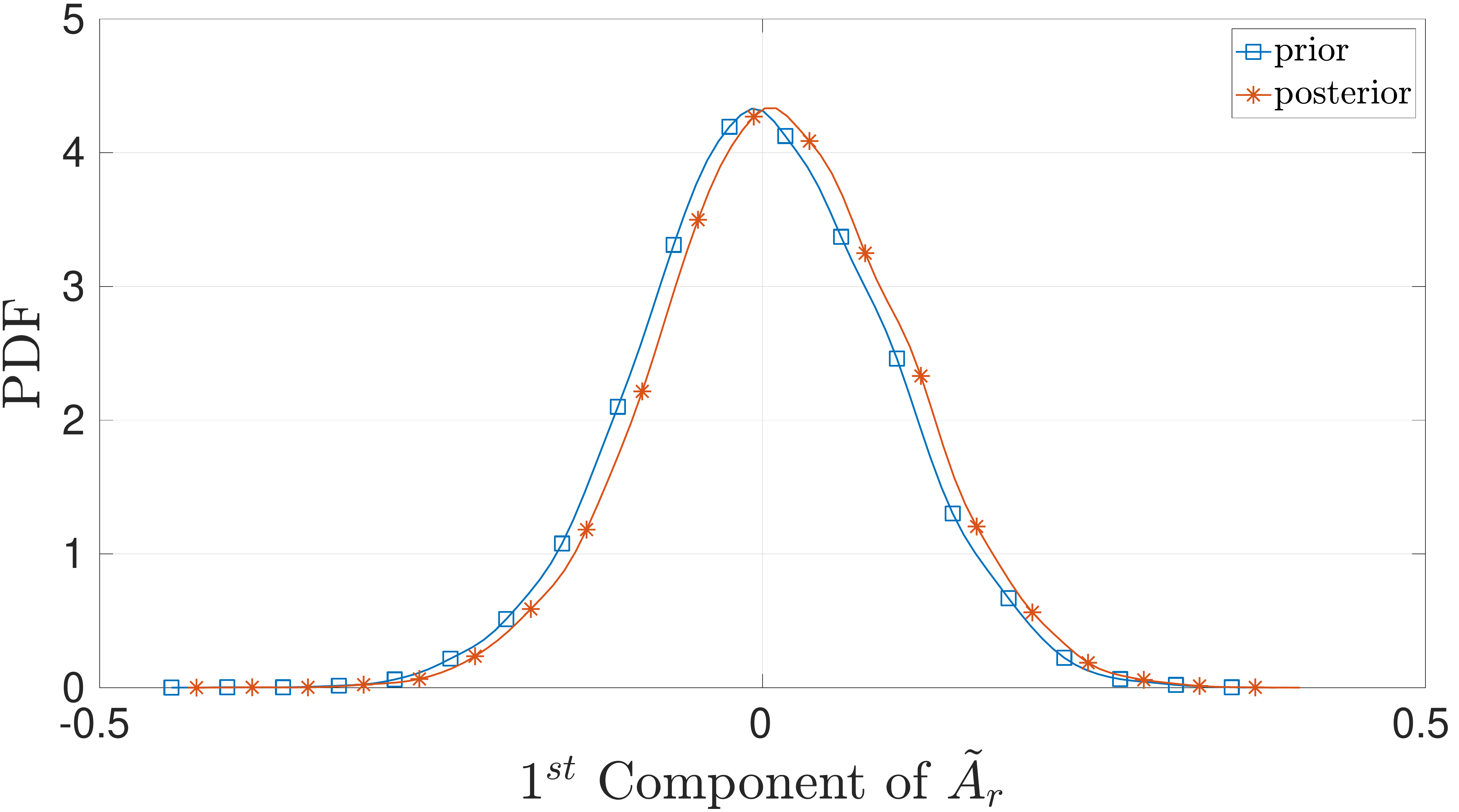}
    \caption{The update of one of the components of $\bar{\bm{A}}_r$}
    \label{fig:Ar_update}
\end{figure}

\begin{figure}
    \centering
    \includegraphics[width=0.45 \textwidth]{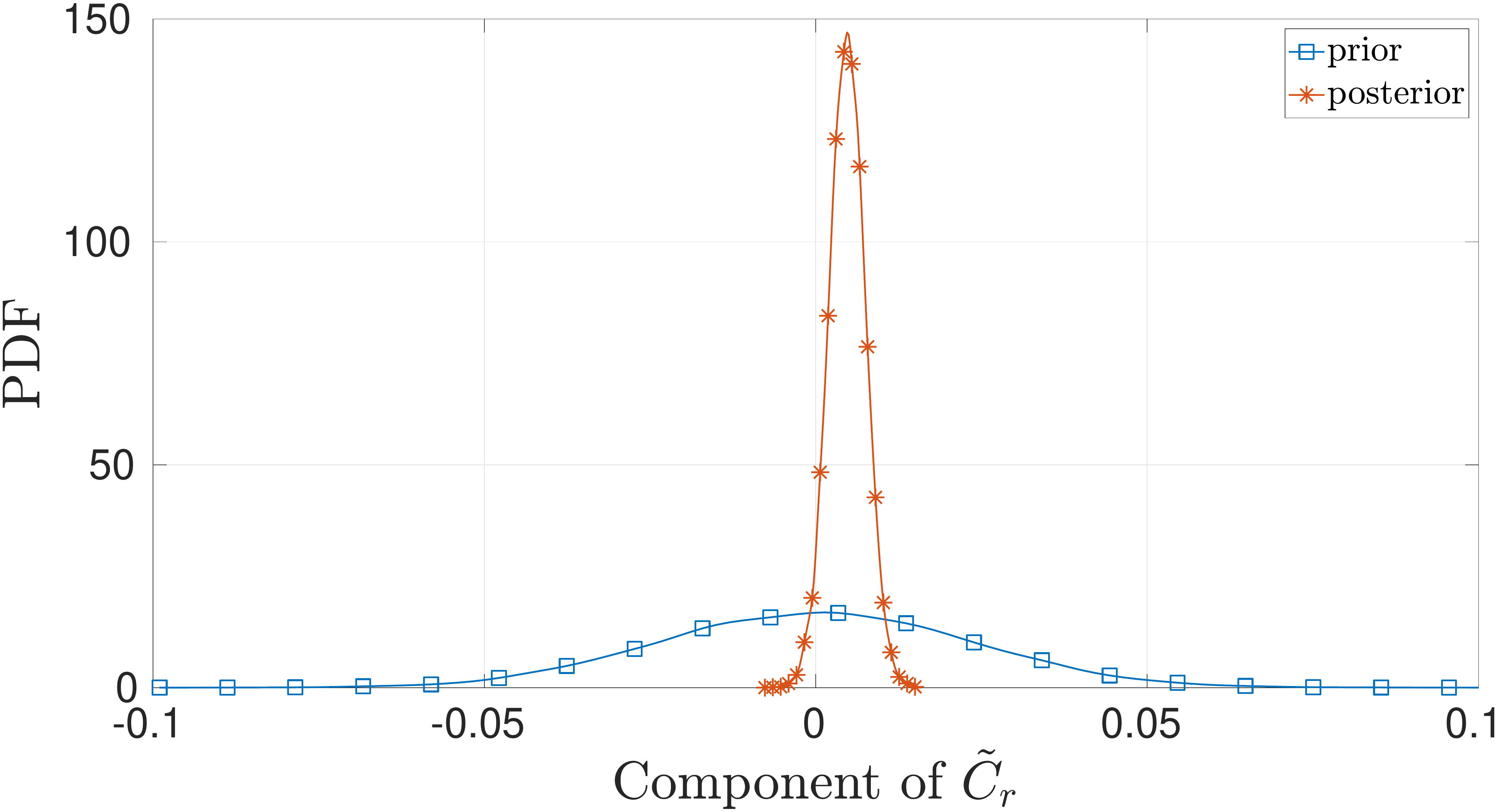}
    \caption{The update of one of the components of $\bar{\bm{C}}_r$}
    \label{fig:Cr_update}
\end{figure}

To make the projection stage independent of the turbulence model, the ROM is first constructed by projecting only the laminar contribution onto the POD basis, see \cite{HiStaMoRo2019,HiAliStaBaRo2018}. Furthermore, the correction terms associated with the turbulence model and possible inaccuracies are identified using the previously described Bayesian approach.
The measurement data are generated by collecting the solution snapshots from $t=60$ seconds to $t=80$ seconds every $\Delta t_{\textrm{train}} = 0.1 \si{s}$, whereas the complete solution is estimated by taking finer time steps, i.e. $\Delta t_s=0.01\si(s)$, in a second order backward differentiation setting. The Navier-Stokes equations are numerically solved using the following spatial discretisation: a Gauss linear scheme for the diffusion terms and a second order upwinding scheme for the convection terms. This setting is used to  store a total number of $200$ snapshots, which are further used to generate the POD basis. In particular the first $10$ modes are used to create a reduced order model according to \autoref{eq:cor_ODE}. For a visual representation of the first $4$ POD modes see \autoref{fig:modes}.
In \autoref{eq:cor_ODE} the reduced operators $\bm A_r$ and $\bm C_r$ are obtained using standard POD-Galerkin projection as explained in \autoref{sec:rom} whereas the correction terms $\bm{\tilde{A}}_r$ and $\bm{\tilde{C}}_r$ are taken to be unknown, and hence described as uncertain. A priori they are modelled as 
Gaussian tensor valued random variables with the zero mean and the variance equal to $1\%$
of the average value. Following \autoref{sec:bayesian}, these are further assimilated with the data (the full order model projected onto the POD basis) by using both 
the EnKF and PCE filters. To avoid fitting the data, the identification is performed only on the segment of the complete time interval, i.e.~by training the correction terms on the first $t_r = [0,4]$s\footnote{Note that this is the time from $60$ to $64$s} where $t_r$ is the time index used to denote the reduced order model simulation. The rest of the interval is used for the validation purposes only. This can be seen in \autoref{fig:coeefs} in which the $99\%$ quantiles of the EnKF posterior for the first $6$ amplitude modes are shown. The figure depicts the time evolution of the temporal amplitude modes, see \autoref{eq:cor_ODE}, before and after the EnKF assimilation. As expected, during the training phase the uncertainty about the amplitude modes is growing in time, and is being significantly reduced right after $4$s when the EnKF smoother is applied. From the results, it is clear that the posterior forecast is following the measurement evolution also outside the training interval. Furthermore, its mean value stays very close to the observation, but is not fully matching due to the discretisation errors, e.g. the posterior is approximated by $1$e$3$ ensemble particles.
To test the convergence, the posterior is estimated by using increasing sequence of ensemble sizes, i.e. $N=1$e$2$, $N=1$e$3$ and $N=1$e$4$, see \autoref{fig:L2EnKF_all}. The largest number of samples is used to estimate the posterior using the complete time interval, $[0\textrm{ }20]$s, in order to build the reference solution. In \autoref{fig:L2EnKF_all} are depicted the relative $L^2$ errors in the posterior velocity, i.e.:
\begin{equation}
    \epsilon(t) = \frac{||\bm{u}_{FOM} - \bm{u}_{ROM}||_{L^2(\mathcal{G})}}{||\bm{u}_{FOM}||_{L^2(\mathcal{G})}},
\end{equation}
in which $ \bm{u}_{ROM}$ is forecasted by taking only the mean of the posterior amplitude modes.  In addition, the error is contrasted to the $L^2$ error of the original snapshots projection onto the POD basis (blue line with diamond markers) and the $L^2$ error of the laminar projection only, here entitled as uncorrected relative error (red line with square markers). The analysis confirmed that, as expected, increasing the number of samples the corrected error reduces and gets closer to the $L^2(\mathcal{G})$ projection of the snapshots onto the POD modes. In contrast, the uncorrected error drastically increases in time, and therefore is not stable. Together, the present findings confirm that an observation time window placed only in the first $4$ seconds of the simulation is sufficient to reproduce the results.
The analysis of the posterior distributions of $\bm \bar{A}_r$ and $\bm \bar{C}_r$ shows the evidence that the diffusion error plays no role in the posterior forecast of the amplitude modes. As can be seen in \autoref{fig:fig1imp}-\autoref{fig:fig2imp} and \autoref{fig:Ar_update}-\autoref{fig:Cr_update}, the first $121$ terms that correspond to the diffusion error $\bm \bar{A}_r$ are having the same posterior variance as the prior one. We believe that this behaviour is justified by the relatively small value of the Reynolds number. In such a case the turbulence model introduces a minimum contribution in terms of additional artificial diffusion. Therefore, these terms are further excluded from the analysis. Next to them, the number of random variables ($\mathcal{O}(N_r^3)$) used to describe the third order convective term is also significantly reduced by excluding non-sensitive variables, the improved ratio of which is near $1$, see \autoref{imp_ratio} and \autoref{fig:fig1imp}. In this manner the total number of $1452$ independent random variables initially describing the problem is reduced to only $93$. The reduction did not impair the accuracy of the identified ROM model as shown in \autoref{fig:L2EnKF_sens}. 
In contrast, a proper selection of the subset is beneficial for both the decrease of the EnKF error, as well as for the PCE identification. 

\begin{figure}
    \centering
    \includegraphics[width=0.45 \textwidth]{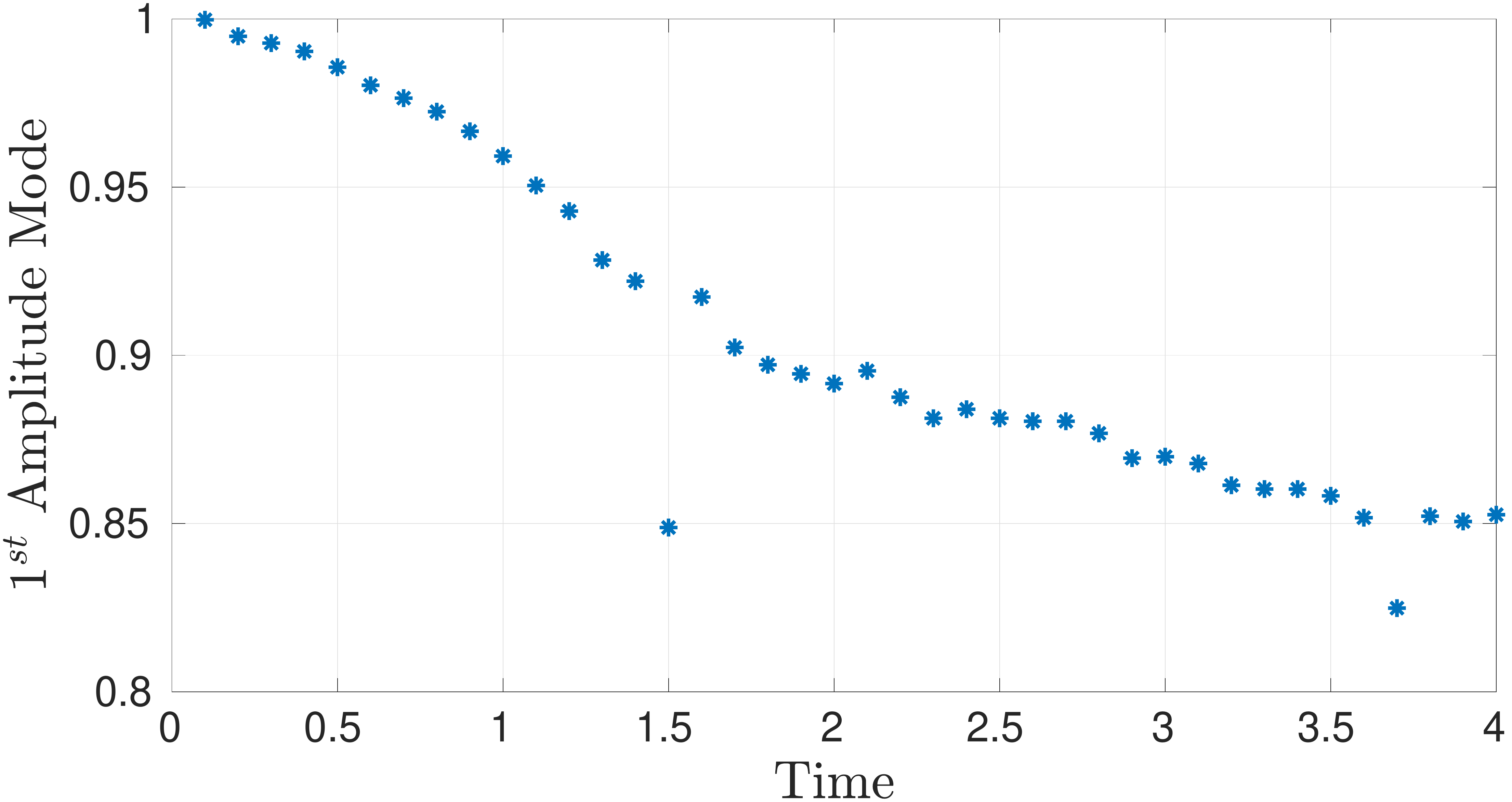}
    \caption{The PCE sparsity of the first amplitude mode}
    \label{fig:sparsity_1st_mode_1e3}
\end{figure}

\begin{figure}
    \centering
    \includegraphics[width=0.45 \textwidth]{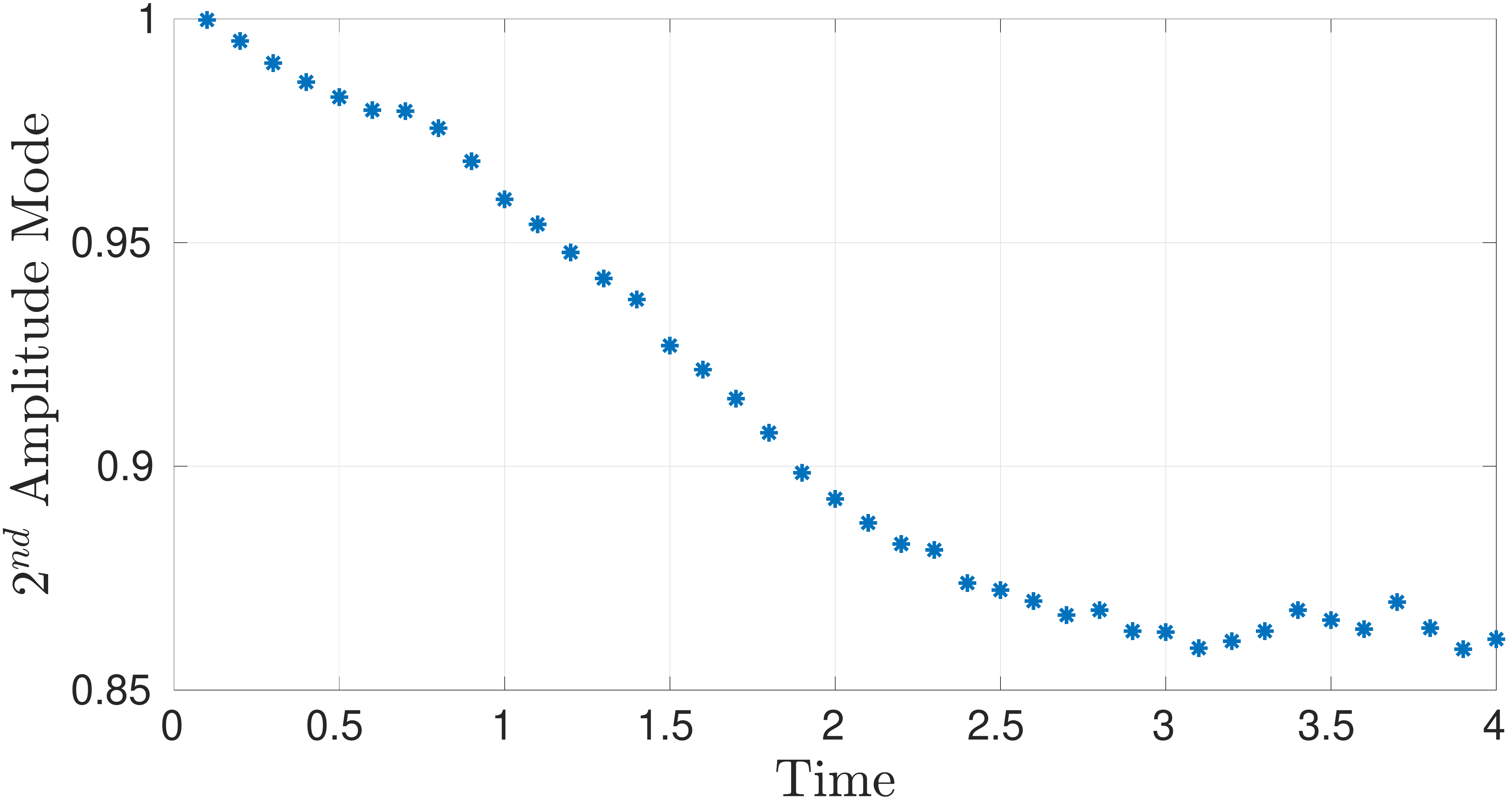}
    \caption{The PCE sparsity of the second amplitude mode}
    \label{fig:sparsity_2nd_mode_1e3}
\end{figure}

\begin{figure}
    \centering
    \includegraphics[width=0.45 \textwidth]{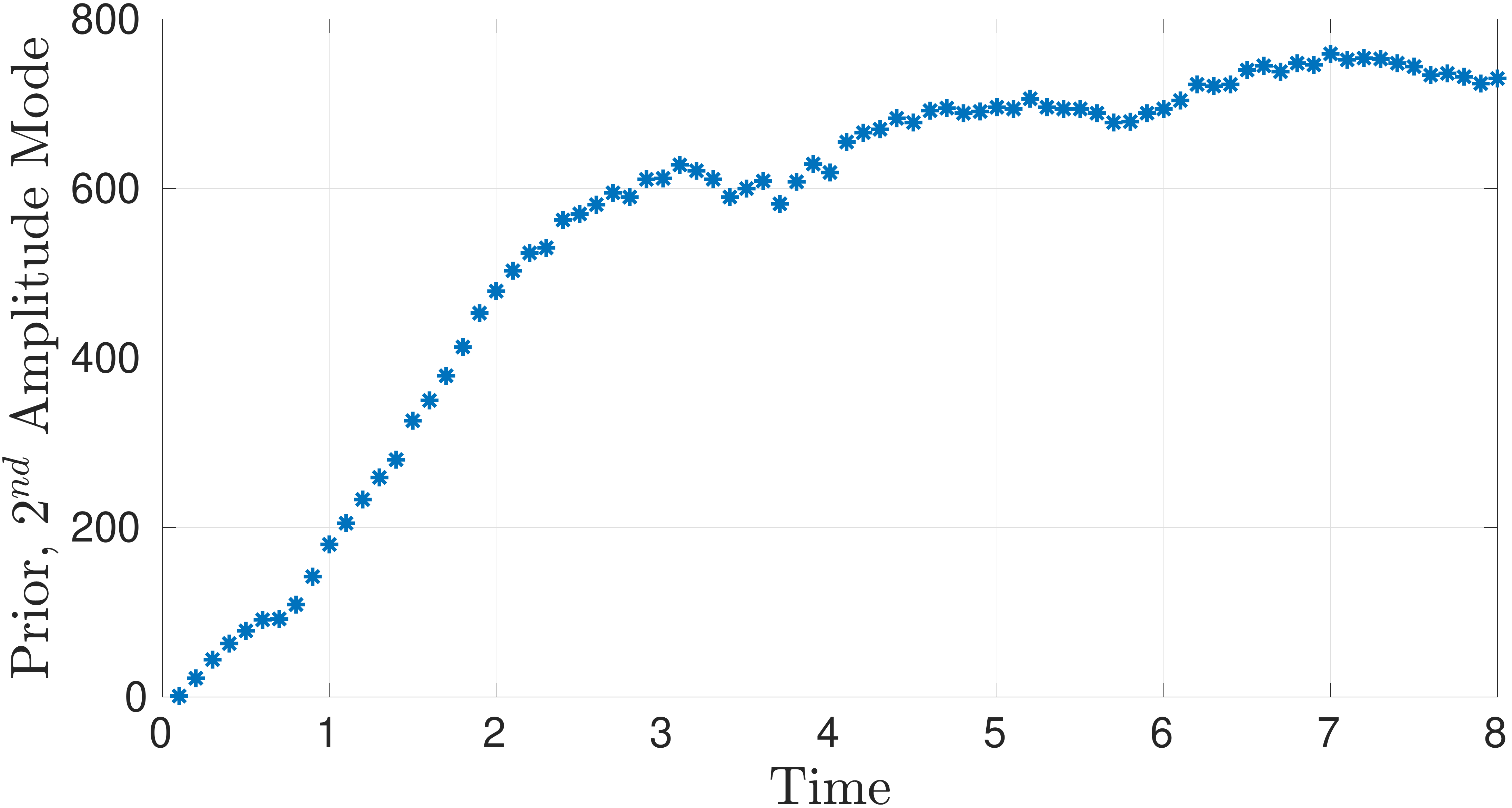}
    \caption{The number of non-zero PCE terms for the second amplitude mode prior to assimilation}
    \label{fig:sparsity_2nd_mode_1e3_num_terms1}
\end{figure}

\begin{figure}
    \centering
    \includegraphics[width=0.45 \textwidth]{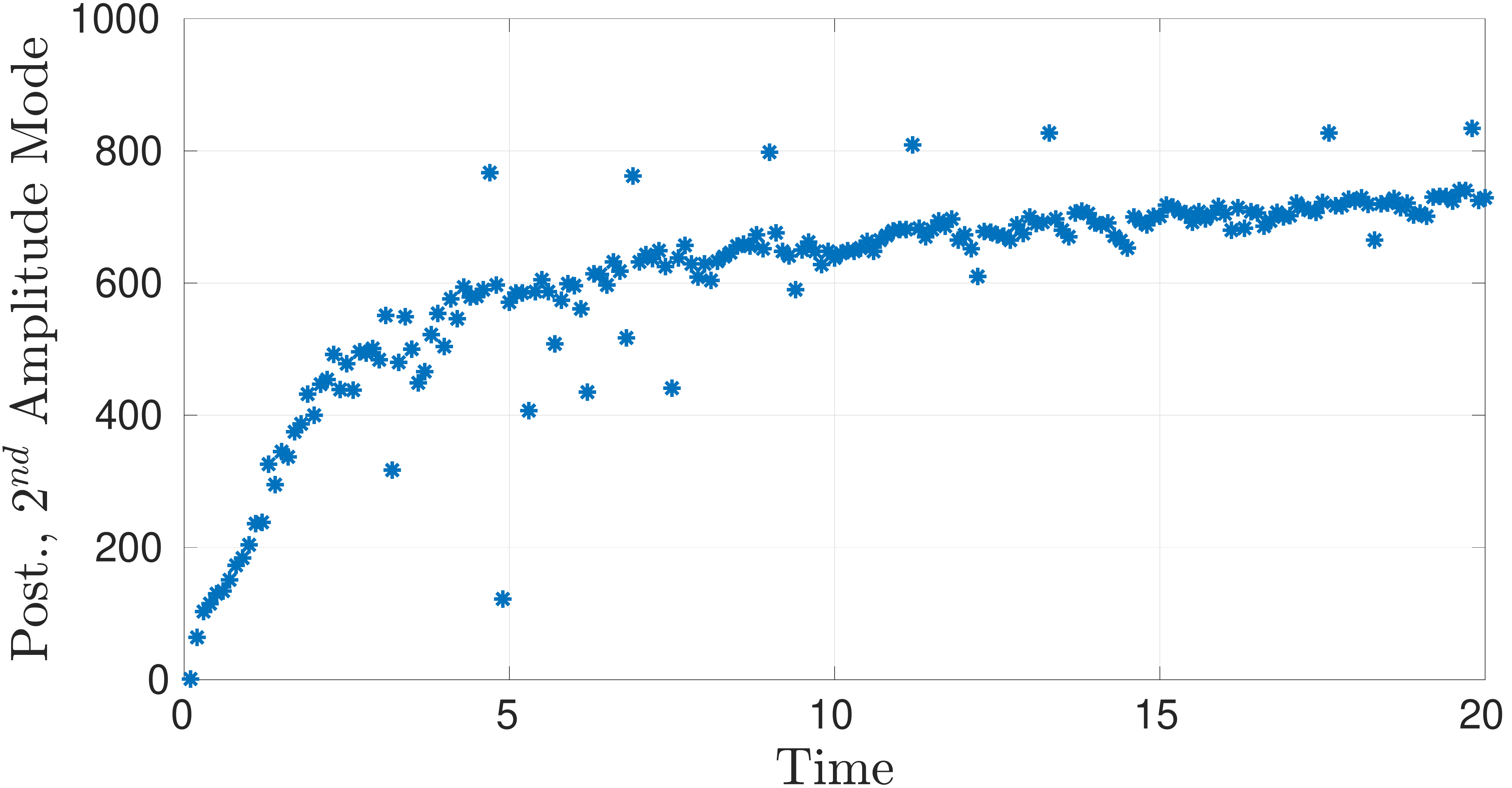}
    \caption{The number of non-zero PCE terms for the second amplitude mode after assimilation}
    \label{fig:sparsity_2nd_mode_1e3_num_terms}
\end{figure}

\begin{figure}
    \centering
    \includegraphics[width=0.45 \textwidth]{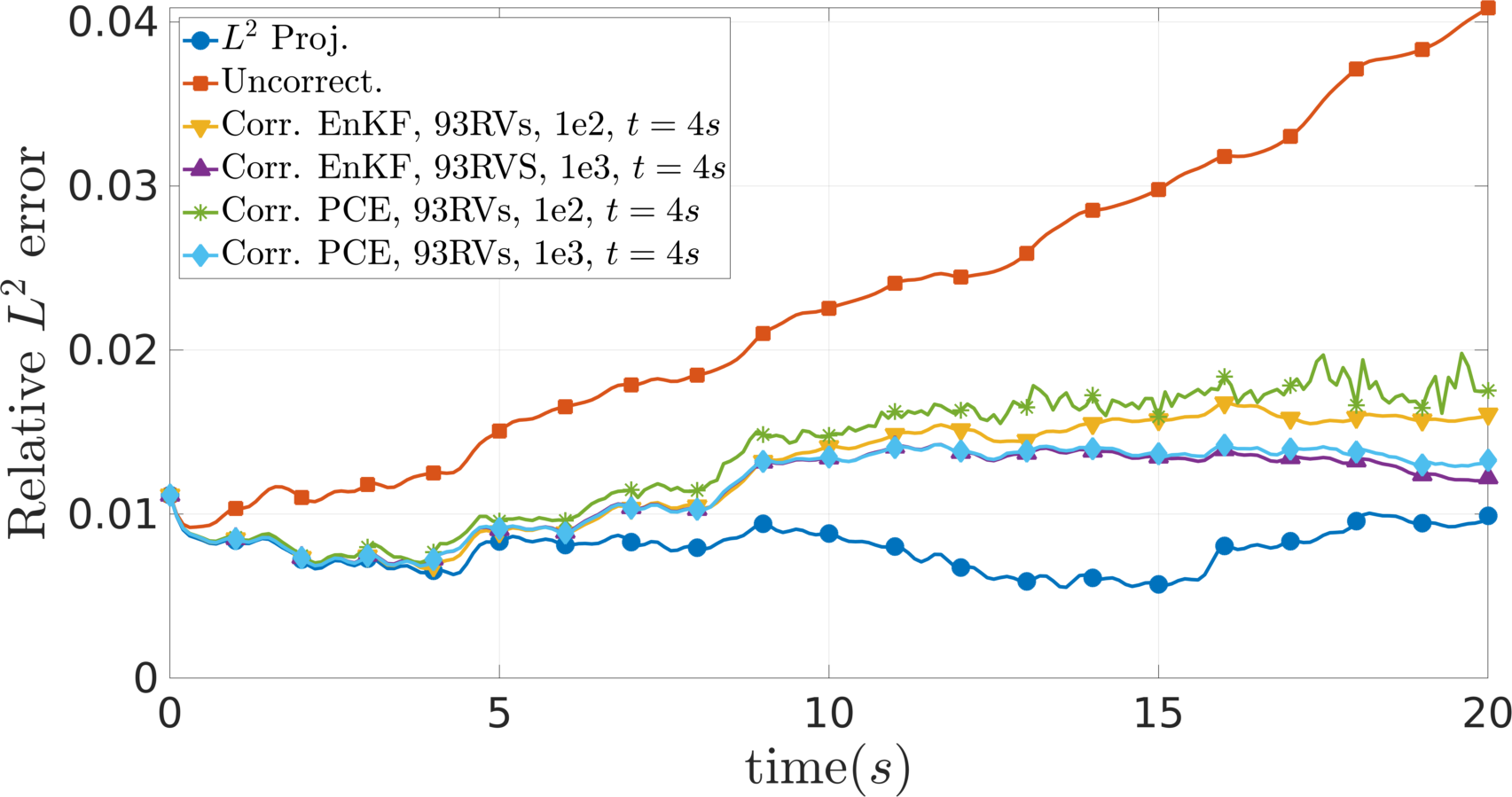}
    \caption{Time evolution of the relative error $\epsilon_u(t)$ of the Bayesian ROM using the EnKF and the PCE techniques on a reduced probability space including $93$ random variables.}
    \label{fig:L2EnKF_PCE}
\end{figure}

The EnKF identification, even though simple, is not robust. The posterior distribution is highly sensitive to the chosen set of ensemble particles. To avoid this, the set of ensembles is replaced by the PCE surrogate model estimated in a non-intrusive way as described in \autoref{subsec:uq}. 
The PCE is built using $93$ random variables of polynomial order $2$ by using $1$e$3$ Monte Carlo samples. This results in $4465$ coefficients describing the posterior. As the posterior cardinality is much higher than the size of the training set, the Bayes's rule with the sparsity prior is employed to compute the coefficients. \autoref{fig:sparsity_1st_mode_1e3}-\autoref{fig:sparsity_2nd_mode_1e3} provide sparsity ratios for the $1$st and $2$nd amplitude mode over the time interval $[0\textrm{ }4]$s. Therefore, one may conclude that only small number of coefficients is different than zero. 
\autoref{fig:L2EnKF_PCE} shows the time evolution of the $L^2$ relative error in velocity for both the EnKF update and the PCE one. From the plots one can deduce that the EnKF and the PCE smoothing, for this particular problem, have similar performances. They give same accuracy over the time interval $[0\textrm{ }10]$s. Afterwords, the accuracy of the PCE slightly deteriorates. The main reason is that the sparsity of the estimate decreases in time, and hence estimation of PCE coefficients requires more samples. In \autoref{fig:sparsity_2nd_mode_1e3_num_terms1}- \autoref{fig:sparsity_2nd_mode_1e3_num_terms} are depicted the sparsities of the $2$nd amplitude mode on the complete time interval before and after assimilation. As analysis shows, the number of non-zero terms grows, and therefore both the number of samples and the polynomial order of the PCE have to be increased. The PCE is known to be sensitive to the long term integration as further explained in \cite{Ro19}. 

However, note that the previous findings were concluded upon one very important assumption. The 
observation is assumed to be subjected to the error 
of Gaussian type $\mathcal{N}(0,\sigma^2)$ with the zero mean and the standard deviation $\sigma$. In a case of the EnKF smoothing, $\sigma$ is taken to be 
$0.1\%$ of the absolute maximum value of the signal. It turns out that smaller values are not possible to be taken as the condition number of the Kalman gain would drastically increase. On the other hand, the PCE smoothing is less prone to such instability. The reason is that the variance is not underestimated, and hence the inverse of the forecasted matrix is well posed. In addition, the standard deviation $\sigma$ following \autoref{bayessparse} can be estimated in a similar manner as PCE coefficients. 

\begin{figure*}
\centering
\begin{minipage}[c]{0.02\textwidth}
\small
$5\si{s}$
\end{minipage}
\begin{minipage}[c]{0.19\textwidth}
\includegraphics[width=\textwidth]{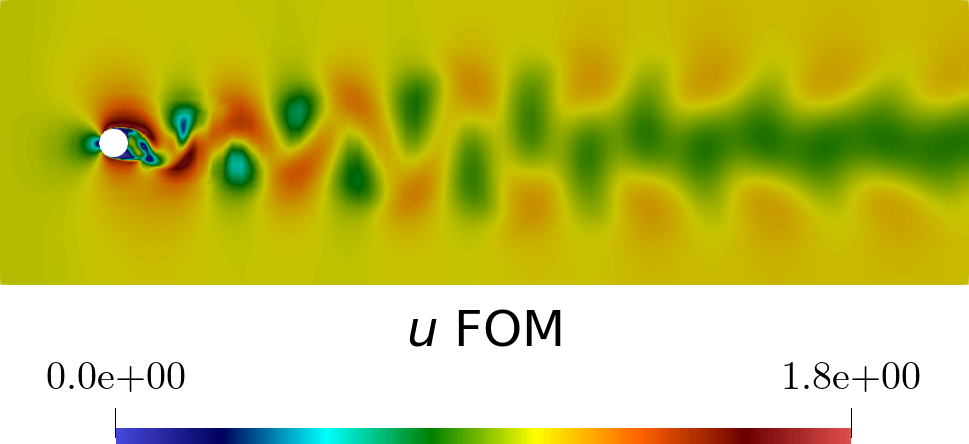}
\end{minipage}
\begin{minipage}[c]{0.19\textwidth}
\includegraphics[width=\textwidth]{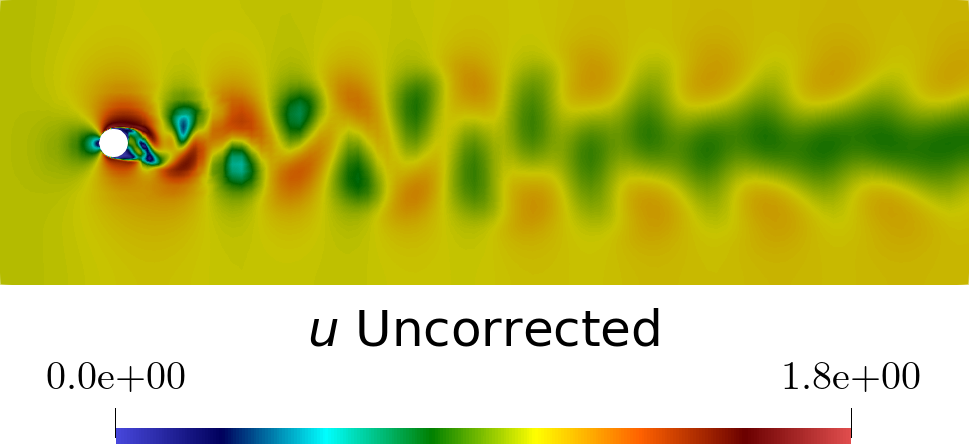}
\end{minipage}
\begin{minipage}[c]{0.19\textwidth}
\includegraphics[width=\textwidth]{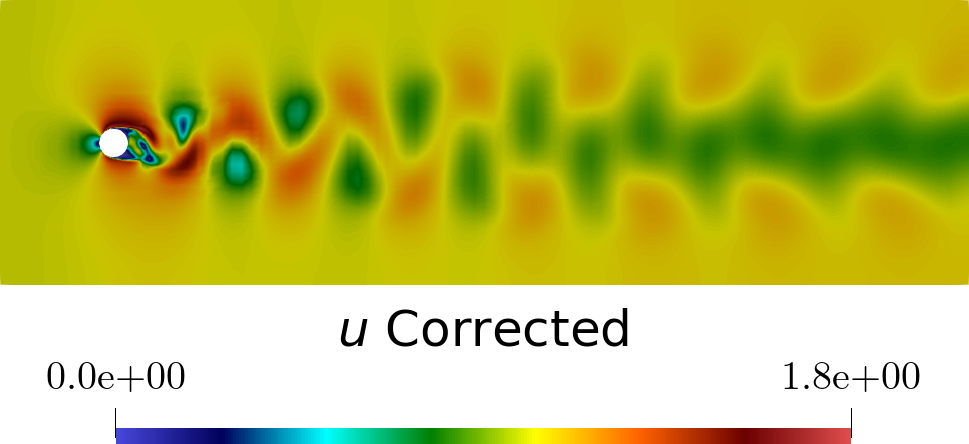}
\end{minipage}
\begin{minipage}[c]{0.19\textwidth}
\includegraphics[width=\textwidth]{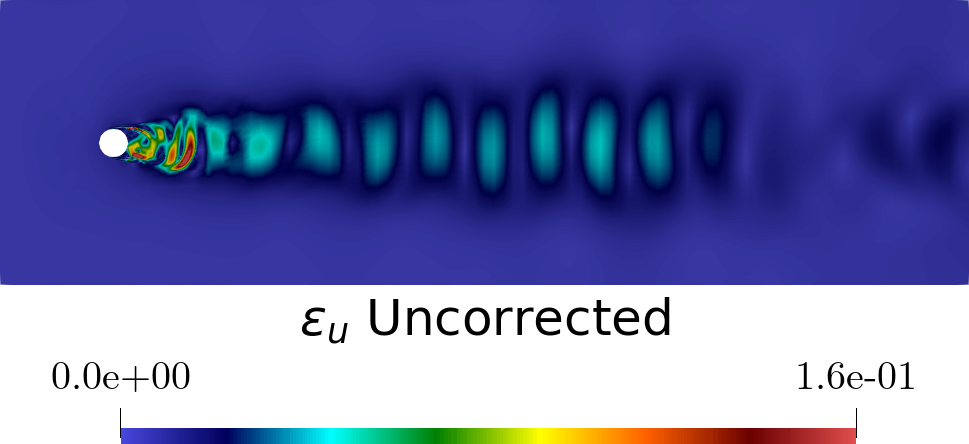}
\end{minipage}
\begin{minipage}[c]{0.19\textwidth}
\includegraphics[width=\textwidth]{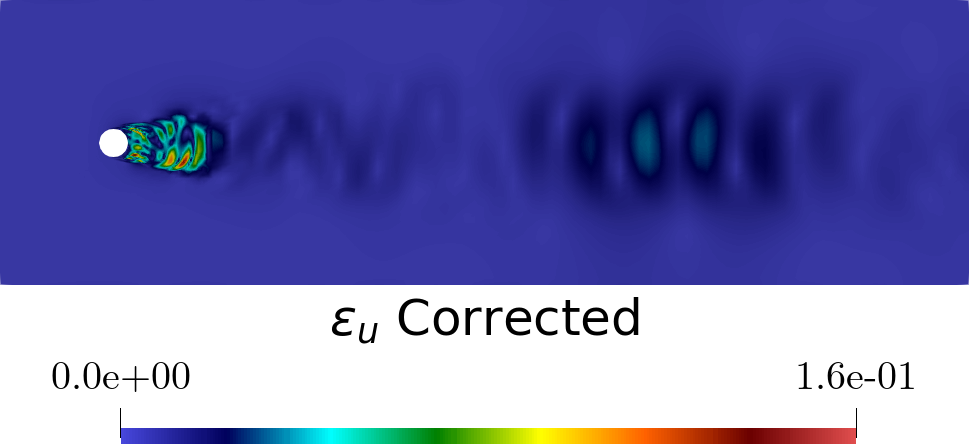}
\end{minipage}
\begin{minipage}[c]{0.02\textwidth}
\small
$10\si{s}$
\end{minipage}
\begin{minipage}[c]{0.19\textwidth}
\includegraphics[width=\textwidth]{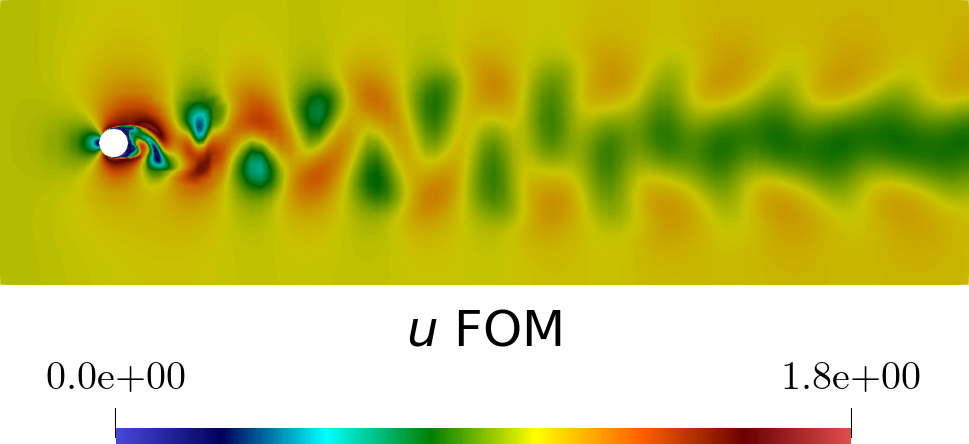}
\end{minipage}
\begin{minipage}[c]{0.19\textwidth}
\includegraphics[width=\textwidth]{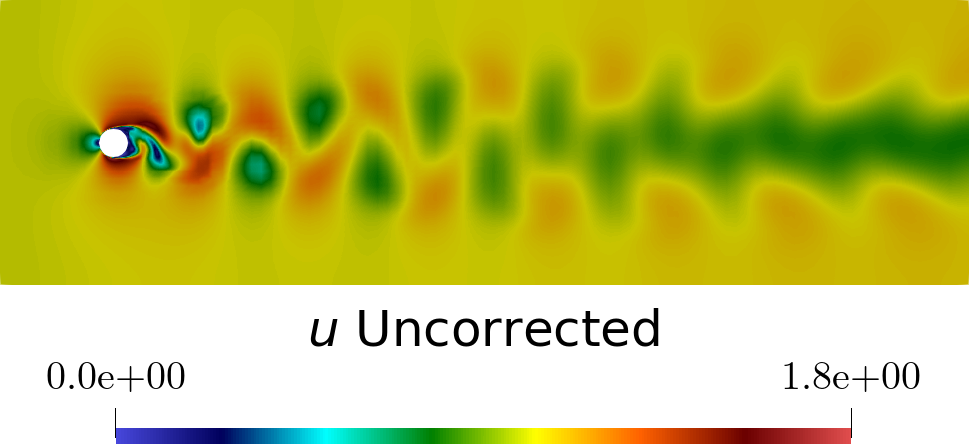}
\end{minipage}
\begin{minipage}[c]{0.19\textwidth}
\includegraphics[width=\textwidth]{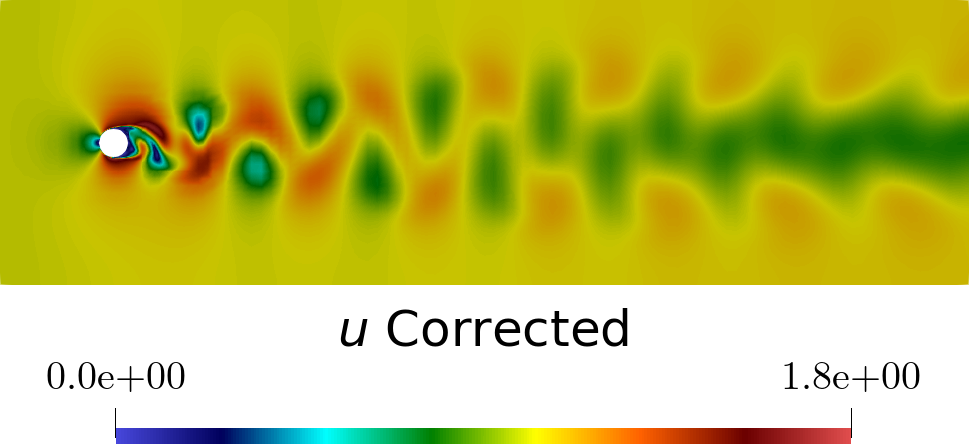}
\end{minipage}
\begin{minipage}[c]{0.19\textwidth}
\includegraphics[width=\textwidth]{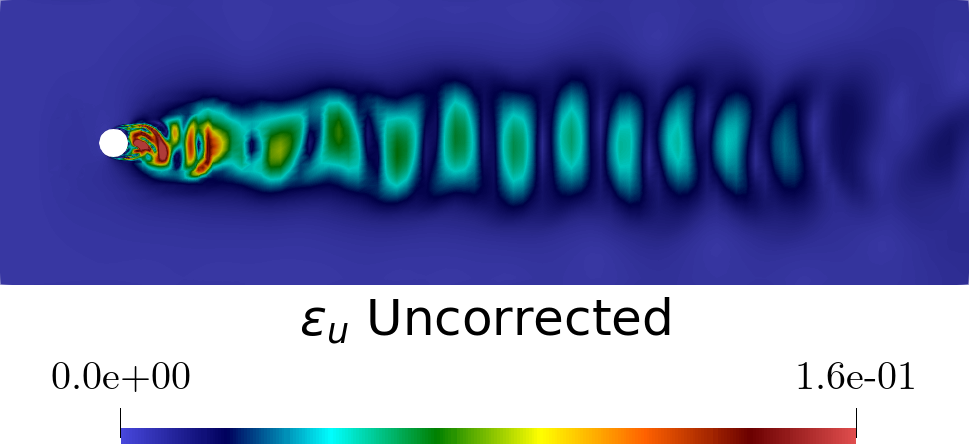}
\end{minipage}
\begin{minipage}[c]{0.19\textwidth}
\includegraphics[width=\textwidth]{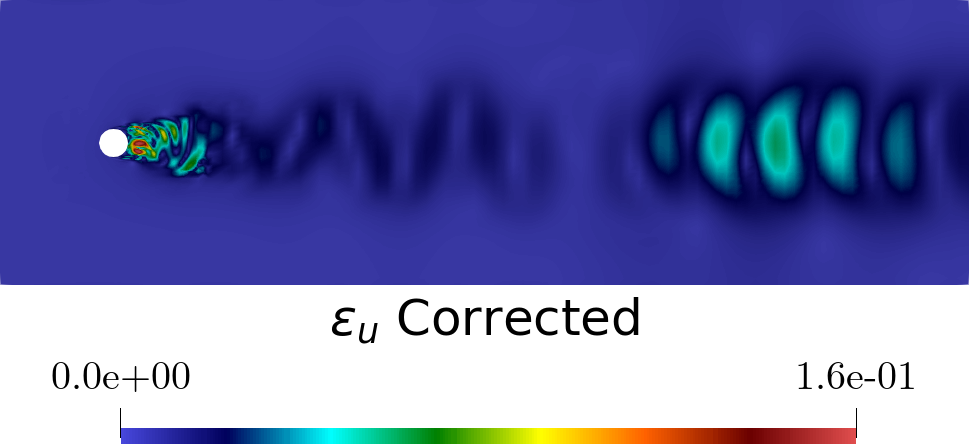}
\end{minipage}
\begin{minipage}[c]{0.02\textwidth}
\small
$15\si{s}$
\end{minipage}
\begin{minipage}[c]{0.19\textwidth}
\includegraphics[width=\textwidth]{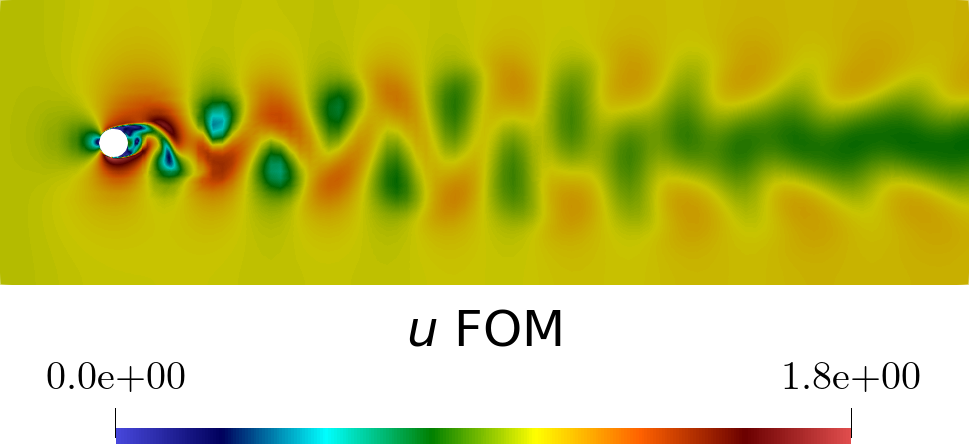}
\end{minipage}
\begin{minipage}[c]{0.19\textwidth}
\includegraphics[width=\textwidth]{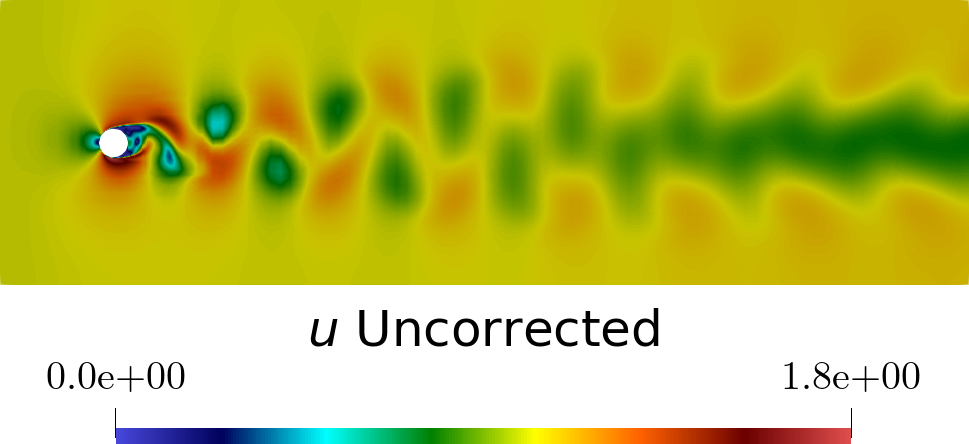}
\end{minipage}
\begin{minipage}[c]{0.19\textwidth}
\includegraphics[width=\textwidth]{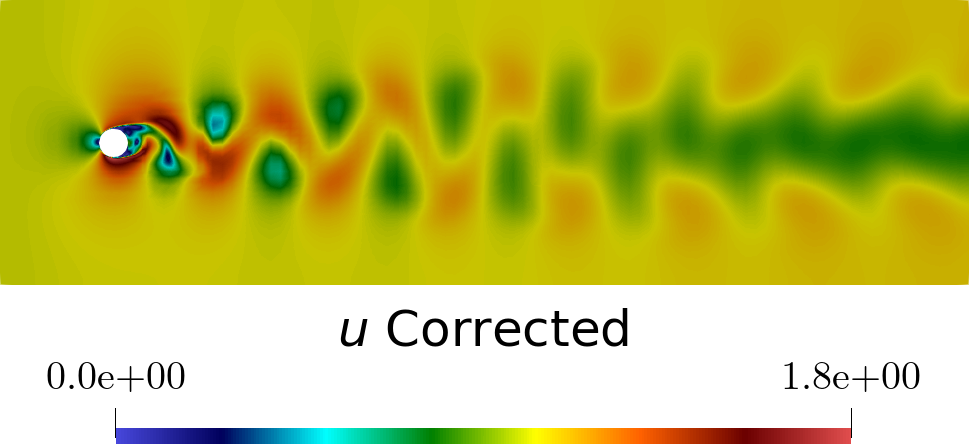}
\end{minipage}
\begin{minipage}[c]{0.19\textwidth}
\includegraphics[width=\textwidth]{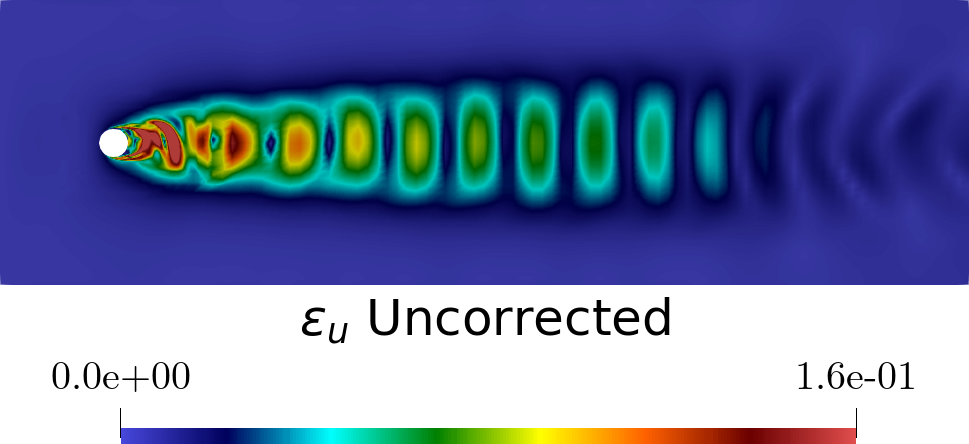}
\end{minipage}
\begin{minipage}[c]{0.19\textwidth}
\includegraphics[width=\textwidth]{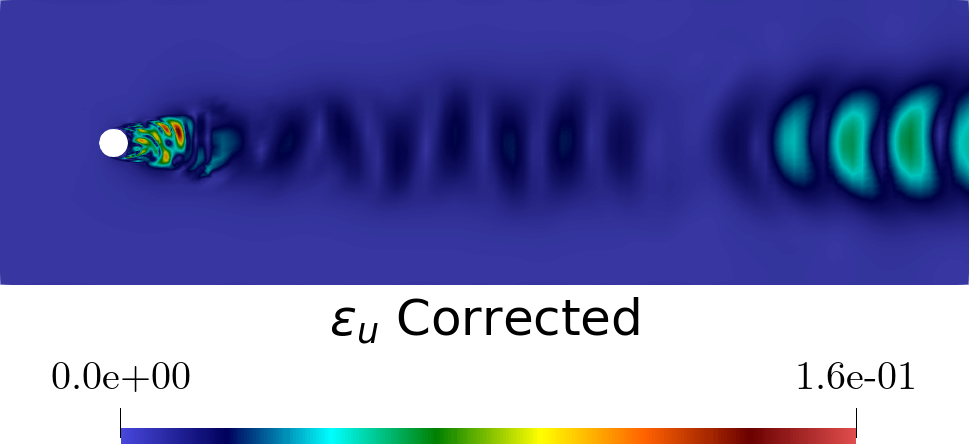}
\end{minipage}
\begin{minipage}[c]{0.02\textwidth}
\small
$20\si{s}$
\end{minipage}
\begin{minipage}[c]{0.19\textwidth}
\includegraphics[width=\textwidth]{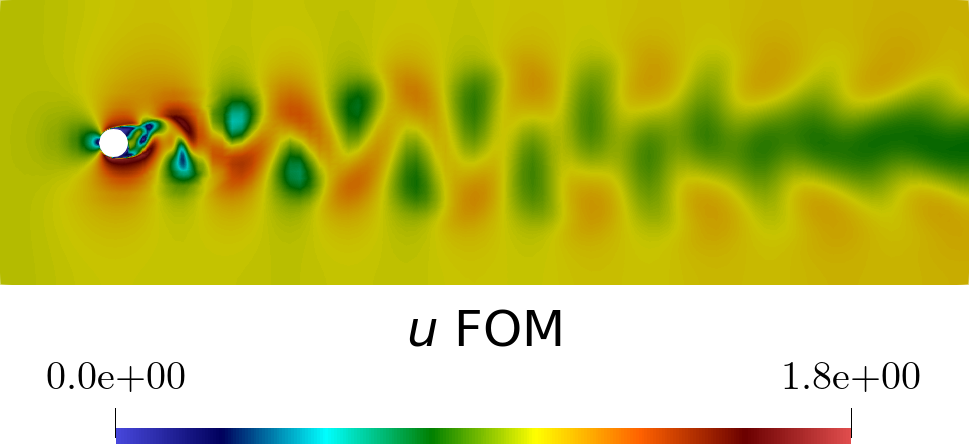}
\end{minipage}
\begin{minipage}[c]{0.19\textwidth}
\includegraphics[width=\textwidth]{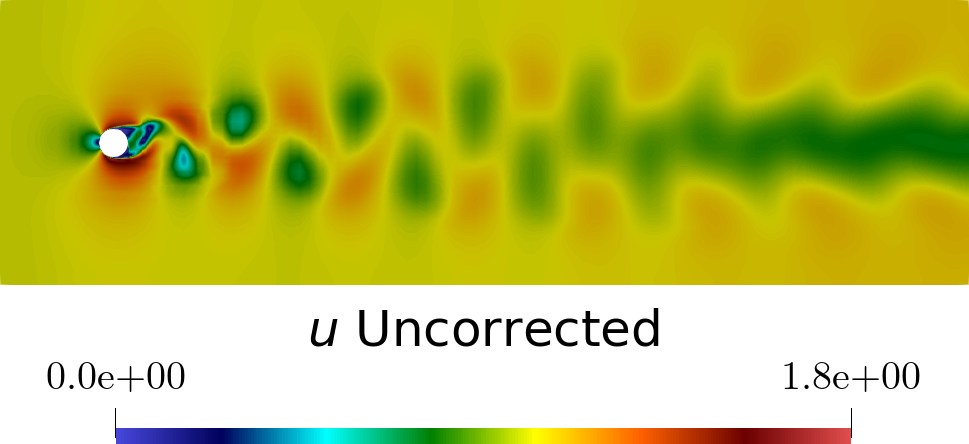}
\end{minipage}
\begin{minipage}[c]{0.19\textwidth}
\includegraphics[width=\textwidth]{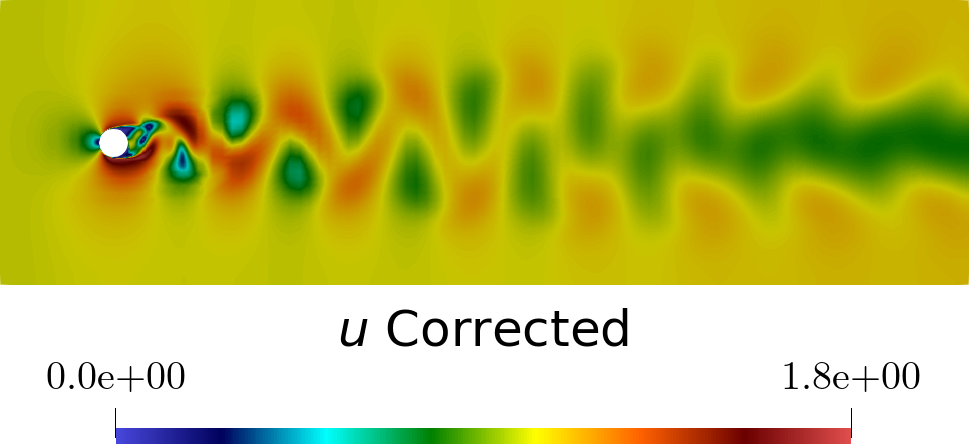}
\end{minipage}
\begin{minipage}[c]{0.19\textwidth}
\includegraphics[width=\textwidth]{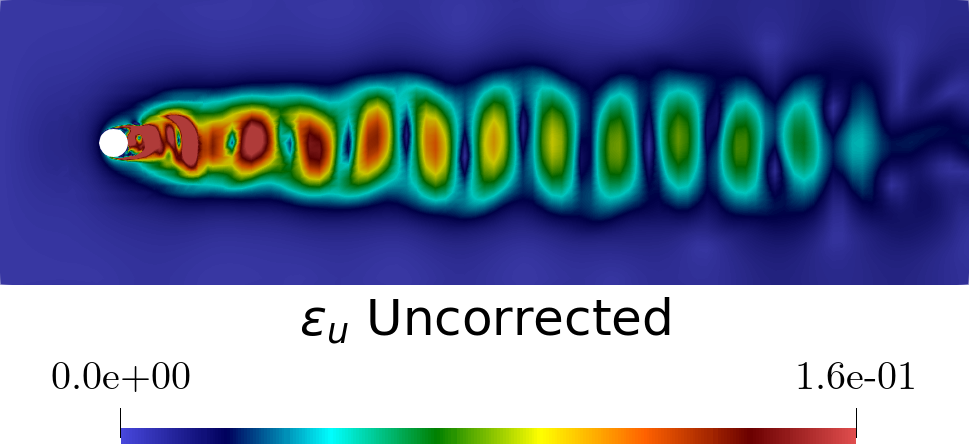}
\end{minipage}
\begin{minipage}[c]{0.19\textwidth}
\includegraphics[width=\textwidth]{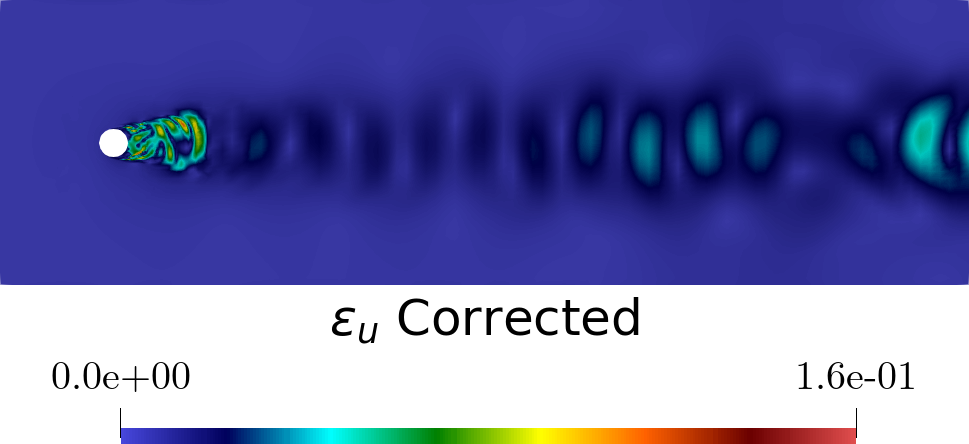}
\end{minipage}
\caption{From the left to right: the FOM-, the uncorrected-, and the  corrected-velocity field using the PCE smoothing with $93$ random variables and $10^3$ samples, the error (w.r.t. FOM) of the uncorrected velocity and the error (w.r.t. FOM) of the corrected velocity field. The different rows, from the top to the bottom, depict different time instants $5\si{s}$,$10\si{s}$,$15\si{s}$ and $20\si{s}$}\label{fig:fields}
\end{figure*}

For a qualitative comparison between the FOM-, the corrected ROM- and the uncorrected velocity field, its several time instances are plotted in \autoref{fig:fields}. Next to these, the corresponding absolute errors are also shown. The results display clear increase of the accuracy of ROM when using the proposed methodology. In contrast to 
the uncorrected ROM, the corrected ROM shows stability with respect to the long time integration.The error is substantially minimised as can be seen in \autoref{fig:L2EnKF_all}, \autoref{fig:L2EnKF_sens} and \autoref{fig:L2EnKF_PCE}.



%% file: tex_files/06_conclusions.tex
\section{Conclusions and future perspectives}\label{sec:conclusions}
In this work a novel approach for the identification of correction terms in a projection-based reduced order model for computational fluid dynamic problems is introduced. The idea of the present approach is to build 
the ROM in which the laminar flow is directly projected onto the POD basis in a Galerkin manner, whereas the turbulence model is incorporated by extending the ROM with the correction terms that are identified given the data in a Bayesian manner. 
 Such an approach permits to obtain not only a deterministic value of the correction terms but also confidence on the estimate. The latter one gives us information about the importance of each of the correction terms, and
 allows further reduction of the model in a probabilistic space. 

The identification is performed using both an EnKF and a PCE smoothing techniques, both of which have shown to reduce the ROM approximation error and to converge properly. In contrast to the direct projection of the snapshots onto the POD basis, the $L_2$ error of the probabilistic approach is more stable and does not increase drastically in time. This is especially true for the PCE model when the appropriate orders are used to approximate the posterior. However, we acknowledge that the PCE estimate is known to fail in long term integration unless used in dynamic adaptive setting, whereas the EnKF is sensitive on the modelling error. From this perspective, both of methods have to be further improved. This is very much the key component in future attempts to deal with the problems that are of
time- and parameter-dependent nature such as the one presented in \cite{GeoStaRoBlu2020}. Looking forward, further attempts could be addressing the accuracy of the POD modes, and their extension to a nonlinear setting.